\documentclass[11pt,a4paper,twoside]{article}
\usepackage{enumerate}
\usepackage{amsmath}
\usepackage{fix-cm}
\usepackage{pb-diagram}
\usepackage[english]{babel}
\usepackage{amssymb,latexsym}
\usepackage{hyperref}
\usepackage[dvips]{color}
\usepackage{authblk}
\def \lra{\longrightarrow}
\def \bui#1#2{\mathrel{\mathop{\kern 0pt#1}\limits^{#2}}}
\def \buil#1#2{\mathrel{\mathop{\kern 0pt#1}\limits_{#2}}}

\let\<\langle
\let\>\rangle
\newcommand{\R}{{\mathbb R}}

\newtheorem{example}{Examples}[section]
\newtheorem{thm}{Theorem}[section]
\newtheorem{lemma}[thm]{Lemma}
\newtheorem{prop}[thm]{Proposition}

\newtheorem{remark}[thm]{Remark}
\newtheorem{remarks}[thm]{Remarks}
\newtheorem{definition}[thm]{Definition}
\newtheorem{notation}[thm]{Notation}
\newtheorem{exabout:ample}[thm]{Example}

\title{Biharmonic Steklov operator on differential forms}

\author[1]{Fida El Chami\thanks{\texttt{fchami@ul.edu.lb}}}
\author[2]{Nicolas Ginoux\thanks{\texttt{nicolas.ginoux@univ-lorraine.fr}}}
\author[1,2]{Georges Habib\thanks{\texttt{ghabib@ul.edu.lb}}}
\author[1]{Ola Makhoul\thanks{\texttt{ola.makhoul@ul.edu.lb}}}

\affil[1]{\footnotesize Lebanese University, Faculty of Sciences II, Department of Mathematics, P.O. Box 90656 Fanar-Matn, Lebanon}
\affil[2]{\footnotesize Universit\'e de Lorraine, CNRS, IECL, F-57000 Metz}

\begin{document}
\maketitle
\noindent\begin{center}\begin{tabular}{p{115mm}}
\begin{small}{\bf Abstract.} 
We introduce the biharmonic Steklov problem on differential forms by considering suitable boundary conditions.
We characterize its smallest eigenvalue and prove elementary properties of the spectrum. 
We obtain various estimates for the first eigenvalue, some of which involve eigenvalues of other problems such as the Dirichlet, Neumann, Robin and Steklov ones.
Independently, new inequalities relating the eigenvalues of the latter problems are proved.
\end{small}\\
\end{tabular}\end{center}

\noindent\begin{small}{\it Mathematics Subject Classification} (2010): 53C21, 58J32, 58C40, 58J50
\end{small}

\noindent\begin{small}{\it Keywords}: Manifolds with boundary, Steklov operator, biharmonic boundary value problem, Robin boundary value problem, eigenvalue estimates  	
\end{small}

\section{Introduction}\label{s:intro}
Let $(M^n,g)$ be an $n$-dimensional compact Riemannian manifold with nonempty smooth boundary $\partial M$. Denote by $\nu$ the inward unit vector field normal to the boundary and by $\Delta f:=-\mathrm{tr}(\nabla^2 f)$ the Laplace operator
applied to a smooth function $f$ on $M$. The following fourth order eigenvalue boundary problem 
\begin{equation}\label{biharmonic-Steklovbis}
\left\{\begin{array}{lll}\Delta^2 f&=0 &\textrm{ on } M\\
f &=0&\textrm{ on }\partial M\\
\Delta f-q \frac{\partial f}{\partial\nu} &=0&\textrm{ on
}\partial M,\end{array}\right.
\end{equation}
also called biharmonic Steklov problem I, or biharmonic Steklov for simplicity, was first introduced by  Kuttler and Sigillito \cite{KS} and Payne \cite{P}. 
Its physical interpretation in terms of the deformation of an elastic plate under the action of transversal forces can be found in e.g. \cite[p. 316]{FGW05} and \cite[p. 2637]{WangXia09}.
When $M$ is a bounded domain in $\R^n$, the spectrum of this problem has been studied in \cite{FGW05} and proved to be discrete consisting of positive eigenvalues of finite multiplicities (see also \cite{BFG} for the case when the boundary is not smooth). 
Variational characterizations of the first eigenvalue $q_1$ have been also given in \cite{FGW05} (see also \cite{Ku}). 
Moreover, in the case of a ball, the spectrum has been calculated explicitly. 

In the case of a compact $n$-dimensional Riemannian manifold $M$ with smooth boundary $\partial M$, sharp estimates for the first eigenvalue $q_1$ of the biharmonic Steklov operator are given in \cite{WangXia09} and \cite{RaulotSavo2015}. It is shown that 
\begin{equation}\label{sharp-est-q_1} q_1 \le \frac{{\rm Vol}(\partial M)}{{\rm Vol}(M)}.\end{equation} 
In \cite[Thm. 2]{RaulotSavo2015}, Raulot and Savo proved  that, if $M$ is a geodesic ball in a space form, then \eqref{sharp-est-q_1} is an equality (see also \cite[Thm. 1.3]{WangXia09}).
Under some assumptions on the Ricci curvature of $M$ and if the mean curvature of the boundary is bounded below by a positive constant $H_0$, then the first eigenvalue $q_1$ also satisfies 
$$q_1 \ge nH_0.$$ 
Moreover, equality holds if and only if $M$ is isometric to a ball of radius $\frac{1}{H_0}$ in $\R^n$ (see \cite[Thm. 1.2]{WangXia09} and \cite[Thm. 2]{RaulotSavo2015}). 
For other recent results on the biharmonic Steklov eigenvalue problem, we refer to \cite{GGS,HS} and the references listed therein.

On the other hand, recall that a compact Riemannian manifold $M$ with smooth nonempty boundary $\partial M$ is called a \emph{harmonic domain} \cite[p. 893]{RS} if and only if it supports a (necessarily unique) solution $f$ to the Serrin boundary value problem
\begin{equation}\label{eq:scalarSerrin}
\left\{\begin{array}{rrl}\Delta f&=1&\textrm{on }M\\ f&=0&\textrm{on }\partial M\\ \partial_\nu f&=c&
\textrm{on }\partial M\end{array}\right.
\end{equation}
for some constant $c\in\R$. 
From \cite{Serrin71,Weinberger71}, we know that the only harmonic domains in $\R^n$ are the Euclidean balls of radius $nc=\frac{n\mathrm{Vol}(M)}{\mathrm{Vol}(\partial M)}$. 
Independently, it is not difficult to check that a solution to the Serrin problem \eqref{eq:scalarSerrin} is an eigenfunction of \eqref{biharmonic-Steklovbis}. \\\\
The aim of this paper is first to extend the biharmonic Steklov problem \eqref{biharmonic-Steklovbis} to the context of differential forms.
As we mentioned above, there is a relationship between problems \eqref{biharmonic-Steklovbis} and \eqref{eq:scalarSerrin}, hence the idea is to also define the Serrin problem on differential forms. 
For this purpose, we assume that
the manifold $M$ carries a non-trivial parallel form and introduce the generalization of \eqref{eq:scalarSerrin} (see \eqref{eq:Serrindiffformsbis}). In this case, we show that $M$ is a harmonic domain if and only if there exists on $M$ a solution to \eqref{eq:Serrindiffformsbis}. 
Also in Section \ref{s:bs}, we provide a natural extension of problem \eqref{biharmonic-Steklovbis} to the case of differential forms (problem \eqref{biharmonic-Steklov}). 
Applying the techniques used in \cite{FGW05} by Ferrero, Gazzola and Weth, we show in Theorem \ref{thm:bihaste} that problem \eqref{biharmonic-Steklov} has a discrete spectrum consisting of a countable number of positive eigenvalues of finite multiplicities. 
This involves proving the ellipticity of problem \eqref{eq:bvporder4} (see Lemma $6.1$ in the appendix). 
Moreover, we give two variational characterizations of the first eigenvalue of this problem (Theorem \ref{vari_BS}) which will be useful to establish inequalities in the sequel. 

In Section \ref{s:ebs}, we obtain different estimates regarding the eigenvalues of problem \eqref{biharmonic-Steklov}. 
As a preliminary step, we prove an interesting property of that problem, namely its invariance by the Hodge star operator. 
On the other hand, under curvature assumptions, using a previous result by Raulot and Savo \cite[Thm. 10]{RaulotSavo2015}, we derive a lower bound \eqref{eq:minq1Theta} for the first eigenvalue $q_{1,p}$, which generalizes the estimates \cite[Thm. 1.2]{WangXia09} and \cite[Thm. 2]{RaulotSavo2015}. 
On the other hand, when the manifold $M$ supports a non-trivial parallel $p$-form, we show that 
$$q_{1,p} \le \frac{{\rm Vol}(\partial M)}{{\rm Vol}(M)},$$
with equality if and only if $M$ is a harmonic domain. 
Surprisingly, when $M$ is a domain of $\R^n$, the eigenvalues of the biharmonic Steklov problem on differential forms are the same as those of the scalar problem, without taking into consideration their multiplicities. 
It should be noted that the same type of result is true for the eigenvalues of the Dirichlet problem.
We end Section \ref{s:ebs} with an inequality relating the eigenvalues corresponding to degrees $p-1$, $p$ and $p+1$ on the sphere. That inequality is established by first proving a more general result when $M$ is isometrically immersed in a Euclidean space and by using the variational characterization \eqref{eq:char1} of the first eigenvalue. The computations involved in these results being rather technical, we place the details in the appendix (Lemma \ref{delta} and Proposition \ref{prop:curvature}) to lighten the text.

In Sections \ref{sec:robinproblem} and \ref{sec:RSvsS}, we establish several bounds concerning the eigenvalues of various differential operators, again using variational characterizations. For example, we give estimates involving the eigenvalues of the Robin eigenvalue problem \eqref{Robin-forms} introduced in \cite{EGH} and those of Neumann \eqref{Neumann-forms}, Dirichlet \eqref{Dirichlet-forms}, biharmonic Steklov \eqref{biharmonic-Steklov} eigenvalue problems on differential forms (Theorems \ref{est_Robin_Dir_Neu} and \ref{est_Robin_Dir_BS}). Note that similar results are known in the case of scalar problems (see \cite[Thm. 1.17]{Ko}). 
In Theorem \ref{t:gabRobin}, we also give an estimate, under some curvature condition along the boundary, for the difference between the first eigenvalues of the Robin problem on $q$ and $(q-p)$-forms, for some $p$ and $q$ such that $p \le q$.\\
Other inequalities involving eigenvalues of the Steklov problem on forms (see Section \ref{sec:RSvsS}) are obtained using properties already established for the eigenvalues and eigenforms of problem \eqref{biharmonic-Steklov}.

{\bf Acknowledgment:} We thank Alberto Ferrero for his skilled guidance in the theory of elliptic boundary value problems and for important references.
The first, third and fourth named benefited from a grant of the Lebanese University.
The second named author thanks the Agence Universitaire de 
la Francophonie (AUF) for its financial support and the Lebanese University for its warm welcome during his stay.
The third named author also thanks the Humboldt Foundation for its support.

\section{Biharmonic Steklov operator on differential forms}\label{s:bs}
\subsection{Serrin problem} \label{sub:serrin}
In this subsection, we extend the Serrin problem to differential forms. This extension will motivate us to define the biharmonic Steklov problem on differential forms.

Recall that the Serrin problem is given by the following \cite{Serrin71}: Let $M\subset \R^n$ be a bounded domain  and let $f$ be a solution to the problem
\begin{equation*}
\left\{\begin{array}{lll}\Delta f&=1 &\textrm{ on } M\\
f &=0&\textrm{ on } \partial M.\\
\end{array}\right.
\end{equation*}
If the inner normal derivative of the function $f$ is a constant $c$, then the domain $M$ must be a ball of radius $nc$ and the function $f$ has the form $(n^2c^2-r^2)/2n$.
Here $c$ is equal to $\frac{{\rm Vol}(M)}{{\rm Vol}(\partial M)}$.
The proof that the Euclidean ball is the unique domain in $\R^n$ supporting a solution to the Serrin problem was given in \cite[Thm. 1]{Serrin71} by using the method of moving planes, which is based on Hopf's maximum principle.
In \cite{Weinberger71}, H.F.~Weinberger suggested an elementary proof introducing so-called $P$-functions for the Laplacian. 
Since then, the Serrin problem has been generalized to several contexts and when the ambient space is a simply connected space form \cite{KP,M,CV}.

A natural question to ask in this set-up is whether the Serrin problem can be extended to differential forms on a domain in $\R^n$. 
For this purpose, we fix $p\in \{0,\cdots,n\}$ and consider, on the set of differential $p$-forms $\Omega^p(M)$, a solution to the system

\begin{equation*}
\left\{\begin{array}{lll}\Delta \omega&=\omega_0 &\textrm{ on } M\\
\omega &=0&\textrm{ on } \partial M,\\
\end{array}\right.
\end{equation*}
where $\omega_0$ is a given parallel form on $\R^n$ assumed to be of norm $1$. 
We now set the following question: If the conditions $\nu \lrcorner d\omega=c\iota^*\omega_0$ and $\iota^*(\delta\omega)=-c\nu\lrcorner\omega_0$ are satisfied on $\partial M$ for some constant $c$ and where $\nu$ is the inner unit normal vector field to $\partial M$, can one deduce that the manifold $M$ is a ball of radius $nc$? Here, $\iota:\partial M\to M$ denotes the inclusion map. Notice that for $p=0$, the problem that we propose reduces to the usual one on functions.
The answer of this question is given in Proposition \ref{prop:serrinformes} below.

\begin{prop} \label{prop:serrinformes}
Let $M$ be a compact manifold with smooth nonempty boundary $\partial M$ and carrying a nontrivial parallel $p$-form $\omega_0$.
\begin{enumerate}
\item\label{statement:0SerrinpSerrin} If $M$ is a harmonic domain, then for the solution $f$ to the Serrin problem on $M$ the $p$-form $\omega:=f\cdot\omega_0$ is the unique solution to the boundary value problem
\begin{equation}\label{eq:Serrindiffformsbis}
\left\{\begin{array}{lll}\Delta\omega&=\omega_0&\textrm{on }M\\
\omega_{|_{\partial M}}&=0&\textrm{on }\partial M\\
\nu\lrcorner d\omega&=c\iota^*\omega_0&\textrm{on }\partial M\\
\iota^*(\delta\omega)&=-c\nu\lrcorner\omega_0&\textrm{on }\partial M
\end{array}\right.
\end{equation}
for some constant $c\in\R$.
\item\label{statement:pSerrin0Serrin} Conversely, if \eqref{eq:Serrindiffformsbis} has a solution $\omega\in\Omega^p(M)$, then assuming w.l.o.g. that $|\omega_0|=1$ on $M$, we have that $\omega=f\cdot\omega_0$ where $f$ solves \eqref{eq:scalarSerrin}.
As a consequence, $M$ must be a harmonic domain.
\end{enumerate}
\end{prop}

{\it Proof.} Before proving the proposition, we begin with the following fact.
Given any parallel form $\alpha$ and a smooth function $h$ on $M$, we have that $\Delta(h\alpha)=(\Delta h)\alpha$. To see this, we first have $d(h\alpha)=dh\wedge\alpha$ and $\delta(h\alpha)=-dh\lrcorner\alpha$, as $\alpha$ is parallel. Therefore, if we take $\{e_1,\cdots, e_n\}$ a local orthonormal frame of $TM$, we compute
\begin{eqnarray*}
\Delta(h\alpha)&=&\delta(dh\wedge\alpha)-d(dh\lrcorner\alpha)\\
&=&-\sum_{i=1}^n e_i\lrcorner(\nabla_{e_i}dh\wedge\alpha)-\sum_{i=1}^n e_i\wedge(\nabla_{e_i}dh\lrcorner\alpha)\\
&=&(\Delta h)\alpha+\sum_{i=1}^n \nabla_{e_i} dh \wedge (e_i\lrcorner\alpha)-\sum_{i=1}^n e_i\wedge(\nabla_{e_i}dh\lrcorner\alpha)\\
&=&(\Delta h)\alpha,
\end{eqnarray*}
since $\nabla dh$ is a symmetric $2$-tensor field.\\
As a first consequence, if $f$ solves \eqref{eq:scalarSerrin}, then for any parallel $p$-form $\omega_0$ on $M$ the $p$-form $\omega:=f\cdot\omega_0$ solves $\Delta\omega=(\Delta f)\cdot\omega_0=\omega_0$ on $M$ together with $\omega_{|_{\partial M}}=f_{|_{\partial M}}\cdot\omega_0=0$.
As for the other two boundary conditions, note that, if $\omega_{|_{\partial M}}=0$, then by $\nabla_X\omega=0$ for all $X\in T\partial M$ we have that $d\omega=\nu^\flat\wedge\nabla_\nu\omega$ and $\delta\omega=-\nu\lrcorner\,\nabla_\nu\omega$.
Therefore,
\[\nu\lrcorner\,d\omega=\nu\lrcorner(\nu^\flat\wedge\nabla_\nu\omega)=\nabla_\nu\omega-\nu^\flat\wedge(\nu\lrcorner\,\nabla_\nu\omega)=\iota^*\nabla_\nu\omega=\partial_\nu f\cdot\iota^*\omega_0=c\iota^*\omega_0\]
and
\[\iota^*(\delta\omega)=-\iota^*(\nu\lrcorner\nabla_\nu\omega)=-\nu\lrcorner\nabla_\nu\omega=-\partial_\nu f\cdot \nu\lrcorner\,\omega_0=-c\nu\lrcorner\,\omega_0,\]
so that $\omega$ solves \eqref{eq:Serrindiffformsbis}.
Note that, since the Dirichlet boundary condition $\omega_{|_{\partial M}}=0$ forces $\ker(\Delta)=\{0\}$ \cite[Thm. p. 445]{Anne89}, the $p$-form $\omega$ is necessarily the only solution to \eqref{eq:Serrindiffformsbis}.
This proves \ref{statement:0SerrinpSerrin}.\\
Conversely, let $\omega$ solve \eqref{eq:Serrindiffformsbis} for some nontrivial parallel $p$-form $\omega_0$.
Up to rescaling $\omega_0$ we may assume that $|\omega_0|=1$ on $M$.
We consider the function $f:=\langle \omega , \omega_0 \rangle$ on $M$.
 By the Bochner formula and $\nabla\omega_0=0$, we have that
$$\Delta f = \langle \nabla^* \nabla \omega , \omega_0 \rangle = \langle \Delta \omega , \omega_0 \rangle - \underbrace{\langle \mathfrak{B}^{[p]}\omega , \omega_0 \rangle}_{0} =|\omega_0|^2=1.$$
Here, we use the fact that the Bochner operator $\mathfrak{B}^{[p]}$ is a symmetric tensor. Also, it is immediate to see that $f_{|_{\partial M}}=0$ since $\omega_{|_{\partial M}}=0$. Therefore, we deduce that
$$\Delta(\omega-f\omega_0)=\omega_0-(\Delta f)\omega_0=0,$$
on $M$ and $(\omega-f\omega_0)|_{\partial M}=0$. Hence by triviality of the Dirichlet kernel, we deduce that $\omega=f\cdot\omega_0$ on $M$. In order to finish the proof, we still have to compute the normal derivative of $f$:
\begin{eqnarray*}
\partial_\nu f &=& \langle \nabla_\nu \omega , \omega_0 \rangle \\
&=&  \langle \iota^* \nabla_\nu \omega , \iota^*\omega_0 \rangle +  \langle \nu \lrcorner \nabla_\nu \omega , \nu \lrcorner \omega_0 \rangle \\
&=& \langle \nu \lrcorner d \omega , \iota^*\omega_0 \rangle -  \langle \iota^*(\delta \omega) , \nu \lrcorner \omega_0 \rangle \\
&=& c \langle \iota^* \omega_0, \iota^*\omega_0 \rangle +c \langle \nu \lrcorner \omega_0, \nu \lrcorner \omega_0 \rangle  \\
&=& c |\omega_0|^2=c.
\end{eqnarray*}
Here we used the identities
$$\left\{
\begin{array}{lll}\nu\lrcorner\nabla_\nu\omega &=& \delta^{\partial M}(\iota^*\omega)- \iota^*(\delta\omega)-S^{[p-1]}(\nu\lrcorner\omega)+(n-1)H\nu\lrcorner\omega \\
\iota^*\nabla_\nu\omega &=& d^{\partial M}(\nu \lrcorner\omega) + \nu \lrcorner d\omega +S^{[p]}(\iota^*\omega).
\end{array}\right.
$$
stated in \cite[Lem. 18]{RaulotSavo2011}. Therefore the function $f$ is a solution to the Serrin problem on $M$.
This concludes the proof of \ref{statement:pSerrin0Serrin} and of Proposition \ref{prop:serrinformes}.
\hfill$\square$

\begin{remark}\rm We notice that if we impose that $\iota^*\omega_0$ is nowhere vanishing along the boundary, the last boundary condition in \eqref{eq:Serrindiffformsbis} can be dropped. 
Indeed, with the boundary condition $\nu\lrcorner d\omega=c\iota^*\omega_0$ and the explicit form $\omega=f\cdot\omega_0$, we compute
\[c\iota^*\omega_0=\nu\lrcorner d\omega=\nu\lrcorner
d(f\omega_0)=(\partial_\nu
f)\omega_0-df\wedge(\nu\lrcorner\omega_0)=(\partial_\nu
f)\iota^*\omega_0,\]
from which $\partial_\nu
f=c$ along $\partial M$ follows. However, the condition $\iota^*\omega_0\neq 0$ is not always assured.
\end{remark}

The Serrin problem on functions is closely related to the biharmonic Steklov operator, that is the boundary problem \eqref{biharmonic-Steklovbis}. 
Indeed, as mentioned in the introduction, on a given compact Riemannian manifold $(M^n,g)$ (not necessarily a domain in $\R^n$), any solution to the Serrin problem is a solution to \eqref{biharmonic-Steklovbis} with $q=\frac{1}{c}$. 
Conversely, it was shown in \cite[Thm. 10]{RS} that the first positive eigenvalue $q_1$ of problem \eqref{biharmonic-Steklovbis} is bounded from below by the first eigenvalue of the Dirichlet-to-Neumann operator on $n$-forms (see Section \ref{sec:RSvsS} for the definition) and, when equality occurs, the corresponding eigenfunction $f$ of \eqref{biharmonic-Steklovbis} is a solution to the Serrin problem. Notice here that, by \cite{FGW05}, problem \eqref{biharmonic-Steklovbis} admits a discrete spectrum that consists of a nondecreasing sequence of positive eigenvalues $(q_i)_i$ of finite multiplicities.

In order to have a similar situation on differential forms, we come back to the Serrin problem defined in \eqref{eq:Serrindiffformsbis} for domains in $\R^n$. In fact, one can easily see
that any solution to the Serrin problem \eqref{eq:Serrindiffformsbis} gives rise to a solution to the following boundary problem

\begin{equation}\label{biharmonic-Steklov}
\left\{\begin{array}{lll}\Delta^2 \omega&=0 &\textrm{ on } M\\ \omega
&=0&\textrm{ on }\partial M\\
\nu\lrcorner\Delta\omega+q\iota^*\delta\omega&=0&\textrm{ on }\partial M\\
\iota^*\Delta \omega-q \nu\lrcorner\,d\omega &=0&\textrm{ on
}\partial M,\end{array}\right.
\end{equation}
with $q=\frac{1}{c}$. The equation $\Delta^2\omega=0$ comes from taking the Laplacian of $\Delta\omega=\omega_0$ and using the fact that $\omega_0$ is a parallel form.
Note here that, because of $\omega_{|_{\partial M}}=0$,   the last two boundary conditions in \eqref{biharmonic-Steklov} are actually equivalent to
\[\Delta\omega=q\nabla_\nu\omega\]
along $\partial M$ since $\nu\lrcorner\,d\omega=\iota^*\nabla_\nu\omega$ and $\iota^*\delta\omega=-\nu\lrcorner\nabla_\nu\omega$ as we have seen in the proof of Proposition \ref{prop:serrinformes}.

It is then natural to study problem \eqref{biharmonic-Steklov} for compact Riemannian manifolds with smooth boundary that are not necessarily domains in $\R^n$.

\subsection{Biharmonic Steklov operator}
In this section, we will show that the spectrum of problem \eqref{biharmonic-Steklov} is discrete and it consists entirely of eigenvalues of finite multiplicities. 
We mainly follow \cite{FGW05}.

First, note that on a compact Riemannian manifold $(M^n,g)$ with smooth boundary $\partial M$, we have the following integration by parts, which is valid for any $\omega,\omega'\in\Omega^p(M)$:
\begin{eqnarray}\label{eq:partialintDelta}
\nonumber&&\int_M\Big(\langle\Delta\omega,\omega'\rangle-\langle\omega,
\Delta\omega'\rangle\Big)\,
d\mu_g\\&=&\int_{\partial M}\Big(\langle\nu\lrcorner
d\omega,\iota^*\omega'\rangle-\langle\iota^*\omega,\nu\lrcorner
d\omega'\rangle+\langle\nu\lrcorner\omega,
\iota^*\delta\omega'\rangle-\langle\iota^*\delta\omega,
\nu\lrcorner\omega'\rangle\Big)\,d\mu_g.
\end{eqnarray}
Thus, replacing $\omega$ by $\Delta\omega$ in
\eqref{eq:partialintDelta}, we obtain:
\begin{eqnarray}\label{eq:partialintDeltacons}
\nonumber\int_M\langle\Delta^2\omega,\omega'\rangle\,d\mu_g
	&=&\int_M\langle\Delta \omega,\Delta\omega'\rangle\,d\mu_g\nonumber\\
&\hspace{-2.5cm}+&\hspace{-1.5cm}\int_{\partial M}\Big(\langle\nu\lrcorner
	d\Delta\omega,\iota^*\omega'\rangle-\langle\iota^*\Delta\omega,
	\nu\lrcorner
	d\omega'\rangle+\langle\nu\lrcorner\Delta\omega,
	\iota^*\delta\omega'\rangle-\langle\iota^*\delta\Delta\omega,
	\nu\lrcorner\omega'\rangle\Big)\,d\mu_g.
\end{eqnarray}
The main result of this section is the following:

\begin{thm}\label{thm:bihaste} Let $(M^n,g)$ be a compact Riemannian manifold with smooth boundary $\partial M$ and let $\nu$ be the inward unit  vector field normal to the boundary. 
Then the boundary problem
\begin{equation*}
\left\{\begin{array}{lll}\Delta^2 \omega&=0 &\textrm{ on } M\\ \omega
&=0&\textrm{ on }\partial M\\
\nu\lrcorner\Delta\omega+q\iota^*\delta\omega&=0&\textrm{ on }\partial M\\
\iota^*\Delta \omega-q \nu\lrcorner\,d\omega &=0&\textrm{ on
}\partial M,\end{array}\right.
\end{equation*}
on $p$-forms, has a discrete spectrum consisting of an unbounded monotonously nondecreasing sequence of positive eigenvalues of finite multiplicities $(q_{j,p})_{j\geq1}$.
\end{thm}

{\it Proof.} As in \cite[Eq. (1.7)]{FGW05}, we let
\[Z:=\left\{\omega\in\Omega^p(M)\,|\,\Delta^2\omega=0\;\textrm{on}
\;M\textrm{ and }\omega_ {
|_{\partial M}}=0\right\}.\]
We define the following Hermitian sesquilinear forms on $Z$: for all
$\omega,\omega'\in\Omega^p(M)$,
\[\left(\omega,\omega'\right)_V:=\int_M\langle\Delta\omega,\Delta\omega'\rangle
\,d\mu_g\;\textrm{and}\;\left(\omega,\omega'\right)_W:=\int_{\partial
M}\langle\nu\lrcorner d\omega,\nu\lrcorner d\omega'\rangle\,
d\mu_g+\int_{\partial
M}\langle\iota^*\delta\omega,\iota^*\delta\omega'\rangle\,d\mu_g.\]
We split the proof of Theorem \ref{thm:bihaste} into the following lemmas.
\begin{lemma} For $\omega \in Z$, we have \begin{equation}\label{partintwithDelta2}\int_M|\Delta\omega|^2\,d\mu_g+\int_{\partial
		M}\Big(\langle\nu\lrcorner\Delta\omega,
	\iota^*\delta\omega\rangle-\langle\iota^*\Delta\omega,
	\nu\lrcorner
	d\omega\rangle\Big)\,d\mu_g=0. \end{equation}
The forms $(\cdot\,,\cdot)_V$ and $(\cdot\,,\cdot)_W$ are positive
	definite on $Z$. 
	Moreover,  there exists a positive constant $C$
	such that $\|\cdot\|_W\leq C\cdot\|\cdot\|_V$ on $Z$.  
	As a consequence, if we denote by $V$ (resp. $W$) the completion of $Z$ w.r.t.
	$\|\cdot\|_V$ (resp. $\|\cdot\|_W$) as Hilbert spaces, then there is a natural
	bounded linear map $I_1\colon V\to W$ extending the identity map $\mathrm{Id}_Z$.
\end{lemma}
{\it Proof.} To prove (\ref{partintwithDelta2}) for $\omega \in Z$, we replace $\omega'$ in (\ref{eq:partialintDeltacons}) by $\omega$ and use the fact that $\Delta^2 \omega=0$ and $\iota^*\omega=	\nu\lrcorner\omega=0$. 
Since, if $\Delta\omega=0$
on $M$ and $\omega_{|_{\partial M}}=0$, then $\omega=0$ (see e.g. \cite[Thm. p. 445]{Anne89}, the sesquilinear form $(\cdot\,,\cdot)_V$ is positive definite. 
For $(\cdot\,,\cdot)_W$, positive definiteness is a consequence of (\ref{partintwithDelta2}). In fact,
if $\left(\omega,\omega\right)_W=0$, then $\nu\lrcorner
d\omega=0$ and $\iota^*\delta\omega=0$ on
$\partial M$ and therefore, from equation (\ref{partintwithDelta2}),  $\Delta\omega=0$ on $M$, from which $\omega=0$ on $M$ follows again by \cite[Thm. p. 445]{Anne89} since $\omega_{|_{\partial M}}=0$. 
We now show,  as in \cite[Sec. 2]{FGW05}, the existence of a positive constant $C$
such that $\|\cdot\|_W\leq C\cdot\|\cdot\|_V$ on $Z$.
First, both $\|\cdot\|_V$ and $\|\cdot\|_{H^2(M)}$ are equivalent on $Z$. To see this, we have for any $\omega\in Z$
\[\|\omega\|_V=\|\Delta\omega\|_{L^2(M)}\leq C\cdot\|\omega\|_{H^2(M)}\]
for some constant $C$ depending only on $M$ and $n$.
On the other hand, by the elliptic estimates and using the fact that, given any
$f\in L^2(M)$, there exists a \emph{unique} weak solution $\omega$
to the boundary value problem $\Delta\omega=f$ on $M$ with $\omega_{|_{\partial
M}}=0$, we have, for any $\omega\in Z$ that $\|\omega\|_{L^2(M)}\leq
C\cdot\|\Delta\omega\|_{L^2(M)}$ for some constant $C$, so that
\[\|\omega\|_{H^2(M)}\leq
C\cdot\left(\|\Delta\omega\|_{L^2(M)}+\|\omega\|_{L^2(M)}
\right)\leq C\cdot\|\Delta\omega\|_{L^2(M)}=C\cdot\|\omega\|_V\]
for some positive constant that we also denote by $C$ and which again depends only on
$M$ and $n$, see e.g. \cite[Thm. 4 in Sec. 6.3]{EvansPDE}.
Therefore, both $\|\cdot\|_V$ and $\|\cdot\|_{H^2(M)}$ are equivalent on $Z$.
Finally, using the fact that both $d$ and $\delta$ are first-order linear differential
operators, we estimate, for any $\omega\in Z$,
\begin{eqnarray*}
\|\omega\|_W^2&=&\|\nu\lrcorner d\omega\|_{L^2(\partial
M)}^2+\|\iota^*\delta\omega\|_{L^2(\partial
M)}^2\\
&\leq&C\cdot\|\omega\|_{H^1(\partial
M)}^2\\
&\leq&C\cdot\|\omega\|_{H^2(M)}^2\\
&\leq&C\cdot\|\omega\|_V^2
\end{eqnarray*}
for some positive constant that we also denote by $C$, which again depends only on
$M$ and $n$.
Here we have also used the boundedness of the trace map $T\colon H^2(M)\to
H^1(\partial M)$.\hfill$\square$\\\\
Next we consider the linear operator $K\colon V\to V$ defined by
\[K:=D_V^{-1}\circ{}^t\!I_1\circ D_W\circ I_1,\]
where $D_V\colon V\to V'$ and $D_W\colon W\to W'$ are the natural duality
isomorphisms, i.e. $D_V(\omega):=\left(\cdot,\omega\right)_V$ and
$D_W(\omega):=\left(\cdot,\omega\right)_W$ for every $\omega$ in $V$ resp. $W$.
As usual, ${}^t\!I_1(\theta):=\theta\circ I_1\in V'$ for every $\theta\in W'$.
Actually $K$ can be defined via the identity
\[\left(K\omega,\omega'\right)_V=\left(I_1\omega,
I_1\omega'\right)_W\left(=\left(\omega,\omega'\right)_W\right)\]
for all $\omega,\omega'\in V$.
By definition, the operator $K$ is
self-adjoint and positive semi-definite. We need now to prove the following.
\begin{lemma}
The map	$I_1$, defined in the previous lemma, is compact and injective. Therefore $K$ is also compact and injective.
\end{lemma}
{\it Proof.}
Let $I_3\colon V\to L^2(\partial M)\oplus L^2(\partial M)$ be the
composition of the following linear maps:
\[
\begin{array}{ccccc}
 V&\longrightarrow& H^{\frac{1}{2}}(\partial M)\oplus
H^{\frac{1}{2}}(\partial M)&\longrightarrow &L^2(\partial M)\oplus L^2(\partial M)\\
\omega&\longmapsto &(\nu\lrcorner
d\omega,\iota^*\delta\omega) &\longmapsto &(\nu\lrcorner d\omega,\iota^*\delta\omega).
\end{array}
\]
Note that $I_3$ is well-defined since the trace operator maps $H^1(M)$ into
(and onto) $H^{\frac{1}{2}}(\partial M)$.
Moreover, since the inclusion map $H^{\frac{1}{2}}(\partial
M)\to L^2(\partial M)$ is compact by the
Rellich-Kondrachov theorem, so is $I_3$.
Now $(Z,\left(\cdot\,,\cdot\right)_W)\to L^2(\partial M)\oplus L^2(\partial
M)$, $\omega\mapsto (\nu\lrcorner
d\omega,\iota^*\delta\omega)$, is a linear isometry, therefore it extends to a
linear isometry $I_2\colon W\to L^2(\partial M)\oplus L^2(\partial M)$ with
$I_3=I_2\circ I_1$.
Since $I_3$ is compact and $I_2$ is a linear isometry, $I_1$ must be compact.\\
We now prove that $I_3$ is injective, so that $I_1$ must be injective as well.
First, we show the inclusion
$$V\subset\left\{\omega\in H^2(M)\cap H_0^1(M)\,|\,\Delta^2\omega=0\textrm{
weakly on }M\right\},$$
where the concept of a weak
solution is defined by the following: given $f\in L^2(M)$, a weak solution $\omega$ to
$\Delta^2\omega=f$ on $M$  is a form $\omega\in H^2\cap H_0^1(M):=H^2(M)\cap
H_0^1(M)$ with
\begin{equation}\label{eq:defweaksol}
\left(\Delta\omega,\Delta\omega'\right)_{
L^2(M) } =\left(f , \omega'\right)_ {
L^2(M)}\;\forall\,\omega'\in H_0^2(M),
\end{equation}
where $H_0^2(M):=\left\{\omega\in H^2(M)\,|\,\omega_{|_{\partial M}}=0\textrm{
and }(\nabla\omega)_{|_{\partial M}}=0\right\}$.
Note that the condition $(\nabla\omega)_{|_{\partial M}}=0$ can be replaced by
$\nabla_\nu\omega=0$ along $\partial M$ because of $\omega_{|_{\partial M}}=0$.
Namely $V\subset H^2(M)$ already holds because of the equivalence of
$\|\cdot\|_V$
and $\|\cdot\|_{H^2(M)}$ on $Z$.
Moreover $V\subset H_0^1(M)$ holds as well because of the
continuous inclusion map $H^2(M)\subset H^1(M)$ and the continuity of the trace
operator $H^1(M)\to L^2(\partial M)$.
Thus $V\subset H^2\cap H_0^1(M)$.
Furthermore, if $\omega\in V$ is given, then there exists a sequence
$(\omega_m)_m$ in $Z$ with $\|\omega_m-\omega\|_V\buil{\lra}{m\to\infty}0$.
Because of $V\subset H^2(M)$ and the equivalence of $\|\cdot\|_V$ and
$\|\cdot\|_{H^2(M)}$, the sequence $(\omega_m)_m$ goes to $\omega$ in $H^2(M)$
and hence $(\Delta\omega_m)_m$ goes to $\Delta\omega$ in $L^2(M)$ .
But since, as a consequence of \eqref{eq:partialintDeltacons}, we have, for all
$m\in\mathbb{N}$ and
$\omega'\in\left\{\omega\in\Omega^p(M)\,|\,\omega_{|_{\partial M}}=0\textrm{
and }(\nabla\omega)_{|_{\partial M}}=0\right\}$,
\[0=\left(\Delta\omega_m,\Delta\omega'\right)_{L^2(M)},\]
we can deduce that $\left(\Delta\omega,\Delta\omega'\right)_{L^2(M)}=0$ for all
$\omega'$ as above and therefore for all $\omega'\in H_0^2(M)$.
This shows that $\omega\in H^2\cap H_0^1(M)$ satisfies $\Delta^2\omega=0$ weakly
on $M$ and therefore the inclusion is proved.
Now, we come back to the injectivity of $I_3$. Consider $\omega\in V$ such that $I_3(\omega)=0$, that is $\nu\lrcorner d\omega=0$ and
$\iota^*\delta\omega=0$ in $L^2(\partial M)$. Then both $d\omega$ and
$\delta\omega$ vanish along $\partial M$ because of $\iota^*\omega=0$,
$\nu\lrcorner\omega=0$ and the identities $[d,\iota^*]=0$ and
$\{\delta,\nu\lrcorner\}=0$.
Again, because of $\omega_{|_{\partial M}}=0$, we have $(d\omega)_{|_{\partial
M}}=\nu^\flat\wedge \nabla_\nu\omega$ and $(\delta\omega)_{|_{\partial
M}}=-\nu\lrcorner\nabla_\nu\omega$, so that $\nu^\flat\wedge
\nabla_\nu\omega=0$ and $\nu\lrcorner\nabla_\nu\omega=0$ along $\partial M$,
from which $\nabla_\nu\omega=0$ on $\partial M$ follows.
This shows that $\omega\in H_0^2(M)$. Taking $\omega'=\omega$ in
\eqref{eq:defweaksol}, we deduce that $\Delta\omega=0$ and therefore
$\omega=0$ on $M$.
This shows $I_3$ and hence $I_1$ to be injective. \hfill$\square$\\
We end the proof of Theorem \ref{thm:bihaste}. Since $K$ is compact, self-adjoint and positive definite in the
Hilbert space $V$, there is a countable Hilbert o.n.b.
$(\omega_i)_{i\geq1}$ of $V$ for which a
monotonously nonincreasing positive real sequence $(\mu_i)_{i\geq1}$ going to $0$
exists
such that $K\omega_i=\mu_i\omega_i$ for all $i\geq1$.
We want to show that, for each $i\geq1$, the eigenform $\omega_i$ lies in $Z$
and satisfies $\mu_i\iota^*\Delta\omega_i=\nu\lrcorner d\omega_i$ as well
as $\mu_i\nu\lrcorner\Delta\omega_i=-\iota^*\delta\omega_i$ along
$\partial M$. 
Hence, for each $i$, the form $\omega_i$ becomes a smooth eigenform for problem
\eqref{biharmonic-Steklov} associated with the eigenvalue $q_{i,p}=\frac{1}{\mu_i}$ which is of finite multiplicity, since $\mu_i$ is. \\
For this purpose, fix $i\geq1$.
Since $\omega_i\in V$, we already know that $\Delta^2\omega_i=0$ holds weakly
on $M$ with $\omega_i{}_{|_{\partial M}}=0$.
It remains to show that $\omega_i$ is smooth and satisfies
$\iota^*\Delta\omega_i=\frac{1}{\mu_i}\nu\lrcorner
d\omega_i$ as well as
$\nu\lrcorner\Delta\omega_i=-\frac{1}{\mu_i}\iota^*\delta\omega_i$ along
$\partial M$.
By definition, for every $\omega\in Z$,
\begin{eqnarray}\label{eq:DeltaomegaDeltaomegaibis}
\mu_i\left(\Delta\omega_i,\Delta\omega\right)_{L^2(M)}
&=&\mu_i\left(\omega_i ,
\omega\right)_V\nonumber\\
&=&\left(K\omega_i,
\omega\right)_V\nonumber\\
&=&\left(\omega_i,\omega\right)_W\nonumber\\
&=&\left(\nu\lrcorner
d\omega_i,\nu\lrcorner d\omega\right)_{L^2(\partial
M)}+\left(\iota^*\delta\omega_i,\iota^*\delta\omega\right)_{L^2(\partial
M)}.
\end{eqnarray}
But by \eqref{eq:partialintDeltacons}, we have, still for every $\omega\in Z$,
\begin{eqnarray}\label{eq:DeltaomegaDeltaomegai}
\nonumber\left(\Delta\omega_i,\Delta\omega\right)_{L^2(M)}&=&\left(\underbrace{
\Delta^2\omega_i}_{0},
\omega\right)_{L^2(M)}-\left(\nu\lrcorner
d\Delta\omega_i,\underbrace{\iota^*\omega}_{0}\right)_{L^2(\partial
M)}+\left(\iota^*\Delta\omega_i,\nu\lrcorner
d\omega\right)_{L^2(\partial M)}\\
\nonumber&&-\left(\nu\lrcorner\Delta\omega_i,
\iota^*\delta\omega\right)_{L^2(\partial
M)}+\left(\iota^*\delta(\Delta\omega_i),\underbrace{\nu\lrcorner\omega}_{0}
\right)_{L^2(\partial M)}\\
&=&\left(\iota^*\Delta\omega_i,\nu\lrcorner d\omega\right)_{L^2(\partial
M)}-\left(\nu\lrcorner\Delta\omega_i,
\iota^*\delta\omega\right)_{L^2(\partial
M)}.
\end{eqnarray}
Here we see both $\iota^*\Delta\omega_i$ and $\nu\lrcorner\Delta\omega_i$ as
elements in $H^{-\frac{1}{2}}(\partial M)$.
Comparing \eqref{eq:DeltaomegaDeltaomegaibis} and \eqref{eq:DeltaomegaDeltaomegai}, we deduce that
\[\left(\iota^*\Delta\omega_i-\frac{1}{\mu_i}\nu\lrcorner d\omega_i,
\nu\lrcorner
d\omega\right)_{L^2(\partial
M)}-\left(\nu\lrcorner\Delta\omega_i+\frac{1}{\mu_i}\iota^*\delta\omega_i,
\iota^*\delta\omega\right)_{L^2(\partial
M)}=0
\]
for all $\omega\in Z$.
Note that the map $Z\to \Omega^p(\partial M)\oplus \Omega^{p-1}(\partial M)$,
$\omega\mapsto(\nu\lrcorner d\omega,\iota^*\delta\omega)$ is continuous w.r.t.
$\|\cdot\|_V$ and
$\|\cdot\|_{L^2(\partial M)}$ and is injective since it is the restriction to
$Z$ of the map $I_3$ from above.
Now, Lemma $6.1$ in the appendix shows that this map is onto. Hence, it follows that
$\iota^*\Delta\omega_i-\frac{1}{\mu_i}\nu\lrcorner
d\omega_i=0$ and
$\nu\lrcorner\Delta\omega_i+\frac{1}{\mu_i}\iota^*\delta\omega_i=0$.
Therefore, $\omega_i$ is
an eigenform for
\eqref{biharmonic-Steklov} associated to the eigenvalue
$q_{i,p}=\frac{1}{\mu_i}$.
The smoothness of $\omega_i$ follows from the fact that both
boundary conditions $\iota^*\Delta\omega_i-\frac{1}{\mu_i}\nu\lrcorner
d\omega_i=0$ and
$\nu\lrcorner\Delta\omega_i+\frac{1}{\mu_i}\iota^*\delta\omega_i=0$ together
with $\omega_i{}_{|_{\partial M}}=0$ define elliptic boundary conditions for
$\Delta^2$.  

In order to finish the proof of Theorem \ref{thm:bihaste}, it remains to show that there is a one-to-one correspondence between solutions of
\eqref{biharmonic-Steklov} and eigenforms of $K$.
We have already shown that every eigenform $\omega$ of $K$,
associated to some eigenvalue $\mu>0$, satisfies  \eqref{biharmonic-Steklov}
with $q=\frac{1}{\mu}$.
Conversely, if $q\in\R$ is given for which a nontrivial solution $\omega$ to
\eqref{biharmonic-Steklov} exists, then by  \eqref{eq:partialintDeltacons}, we
have, for every $\omega'\in Z$,
\[\left(\Delta\omega,\Delta\omega'\right)_{L^2(M)}
=q\cdot\left(\left(\nu\lrcorner
d\omega,\nu\lrcorner
d\omega'\right)_{L^2(\partial
M)}+\left(\iota^*\delta\omega,\iota^*\delta\omega'\right)_{L^2(\partial
M)}\right),\]
that is, seeing both $\omega,\omega'$ as elements of $V$,
\[\left(\omega,\omega'\right)_V=q\cdot\left(\omega,\omega'\right)_W.\]
Note that necessarily $q>0$ holds, otherwise $\omega=0$ would follow.
By definition of $K$, we then have
$\left(\omega,\omega'\right)_V=q\cdot\left(K\omega,\omega'\right)_V$ for all
$\omega'\in Z$ and hence in $V$, therefore $K\omega=\frac{1}{q}\omega$.
This shows $\omega$ to be an eigenform of $K$ associated to the eigenvalue
$\mu=\frac{1}{q}$.
This shows the correspondence to be one-to-one. This concludes the proof.
\hfill$\square$

In the following, we give a characterization for the first eigenvalue $q_{1,p}$ on $p$-forms. This will be used later in order to get estimates for the eigenvalues.

\begin{thm}\label{vari_BS} The first eigenvalue $q_{1,p}$ of the boundary problem \eqref{biharmonic-Steklov}
is characterized by

\begin{eqnarray}
\label{eq:char1}q_{1,p}&=&\inf\left\{\frac{\|\Delta\omega\|_{L^2(M)}^2}{
\|\nu\lrcorner d\omega\|_{L^2(\partial
M)}^2+\|\iota^*\delta\omega\|_{L^2(\partial
M)}^2}\,|\,\omega\in\Omega^p(M),\,\omega_{|_{\partial
M}}=0\textrm{ and }\nabla_\nu\omega\neq0\right\}\\
\label{eq:char2}&=&\inf\left\{\frac{\|\omega\|_{L^2(\partial M)}^2}{
\|\omega\|_{L^2(
M)}^2}\,|\,\omega\in\Omega^p(M)\setminus\{0\},\,\;\Delta\omega=0\textrm{ on
}M\right\}.
\end{eqnarray}
Both infima are indeed minima, \eqref{eq:char1} is attained by an eigenform of \eqref{biharmonic-Steklov}, associated to $q_{1,p}$ and \eqref{eq:char2} is attained by $\Delta \omega$, where $\omega$ is an eigenform of \eqref{biharmonic-Steklov}, associated to $q_{1,p}$.
\end{thm}

{\it Proof.} As mentioned above, it follows from \eqref{partintwithDelta2} that,
given any
nonzero eigenform $\omega$ associated to a positive eigenvalue $q$ of
\eqref{biharmonic-Steklov}, we have
\[\int_M|\Delta\omega|^2\,d\mu_g=q\int_{\partial M}|\nu\lrcorner
d\omega|^2\,d\mu_g+q\int_{\partial M}|\iota^*\delta\omega|^2\,d\mu_g,\]
so that $q_{1,p}\leq \frac{\|\Delta\omega\|_{L^2(M)}^2}{
\|\nu\lrcorner d\omega\|_{L^2(\partial
M)}^2+\|\iota^*\delta\omega\|_{L^2(\partial
M)}^2}$ for every such eigenform, with equality for $\omega$ associated to
$q_{1,p}$ of course.
More generally, for every $\omega\in V$, we can write
$\omega=\sum_i\left(\omega,\omega_i\right)_V\cdot\omega_i$ because
$(\omega_i)_i$ is a Hilbert orthonormal basis of $V$.
From this, we can express
\[\|\Delta\omega\|_{L^2(M)}^2=\|\omega\|_V^2=\sum_i\left|\left(\omega ,
\omega_i\right)_V\right|^2\]
on the one hand, and
\begin{eqnarray*}
\|\nu\lrcorner d\omega\|_{L^2(\partial
M)}^2+\|\iota^*\delta\omega\|_{L^2(\partial
M)}^2&=&\|\omega\|_W^2\\
&=&\left(K\omega,
\omega\right)_V\\
&=&\sum_i\frac { 1 }{ q_{ i , p } } \left|
\left(\omega,\omega_i\right)_V\right|^2\\
&\leq&\frac{1}{q_{1,p}}\sum_i\left|
\left(\omega,\omega_i\right)_V\right|^2\\
&\leq&\frac{1}{q_{ 1 ,
p}}\|\Delta\omega\|_{
L^2(M)}^2
\end{eqnarray*}
on the other hand, therefore
\[q_{1,p}=\inf\left\{\frac{\|\Delta\omega\|_{L^2(M)}^2}{
\|\nu\lrcorner d\omega\|_{L^2(\partial
M)}^2+\|\iota^*\delta\omega\|_{L^2(\partial
M)}^2}\,|\,\omega\in V\setminus\{0\}\right\}.\]
Now, we will show that this infimum can be taken over all $\omega\in(H^2\cap
H_0^1)(M)\setminus\{0\}$, as well as over all smooth forms vanishing and whose normal derivative does not vanish identically along
$\partial M$.
Recall that,
\[H_0^2(M)=\left\{\omega\in (H^2\cap
H_0^1)(M)\,|\,\iota^*\delta\omega=0\textrm{ and }\nu\lrcorner
d\omega=0\right\},\]
since, as we noticed above,
if $\omega\in (H^2\cap H_0^1)(M)$ is such that both $\iota^*\delta\omega$ and
$\nu\lrcorner d\omega$ vanish along $\partial M$, then so does
$\nabla_\nu\omega$.
As in \cite[Thm. 1.2]{FGW05}, we have the following lemma:
\begin{lemma} The inner product $\left(\cdot\,,\cdot\right)_V$ is well defined
on $(H^2\cap H_0^1)(M)$ and we have the following
$\left(\cdot\,,\cdot\right)_V$-orthogonal splitting:
\[H^2\cap H_0^1(M)=V\oplus
H_0^2(M).\]
\end{lemma}
{\it Proof.} By its definition, $\left(\cdot\,,\cdot\right)_V$ is well defined
on $(H^2\cap H_0^1)(M)$.
Furthermore, \eqref{eq:partialintDeltacons} already implies
that, for all
$(\omega,\omega')\in Z\times H^2\cap H_0^1(M)$,
\begin{equation}\label{eq:innerprodV}
\left(\omega,
\omega'\right)_V=\left(\iota^*\Delta\omega , \nu\lrcorner
d\omega'\right)_{L^2(\partial M)}-\left(\nu\lrcorner\Delta\omega
,\iota^*\delta\omega'\right)_{L^2(\partial M)},
\end{equation}
so that $\left(\omega,\omega'\right)_V=0$ as soon as $\omega'\in H_0^2(M)$.
This shows that $H_0^2(M)\subset Z^\perp=V^\perp$.
Conversely, let $\omega'\in V^\perp\subset H^2\cap H_0^1(M)$.
Then $\left(\omega_i,\omega'\right)_V=0$ for all $i$.
By \eqref{eq:innerprodV}, this is equivalent to $\left(\nu\lrcorner
d\omega_i,\nu\lrcorner d\omega'\right)_{L^2(\partial
M)}+\left(\iota^*\delta\omega_i,\iota^*\delta\omega'\right)_{L^2(\partial
M)}=0$ for all $i$ by
$\iota^*\Delta\omega_i=\frac{1}{\mu_i}\nu\lrcorner\,d\omega_i$ and
$\nu\lrcorner\Delta\omega_i=-\frac{1}{\mu_i}\iota^*(\delta\omega_i)$.
Since the map $Z\to \Omega^p(\partial M)\oplus\Omega^{p-1}(\partial M)$,
$\omega\mapsto(\nu\lrcorner d\omega,\iota^*\delta\omega)$
is bounded (w.r.t. $\|\cdot\|_V$ and $\|\cdot\|_{L^2(\partial M)}$), onto
and the $\omega_i$'s span a dense subspace of $V$, we obtain that
$\left(\omega_1,\nu\lrcorner d\omega'\right)_{L^2(\partial M)}=0$ as well as
$\left(\omega_2,\nu\lrcorner d\omega'\right)_{L^2(\partial M)}=0$ for all
$(\omega_1,\omega_2)\in \Omega^p(\partial M)\oplus\Omega^{p-1}(\partial M)$ and
therefore $\nu\lrcorner d\omega'=0$ and $\iota^*\delta\omega'=0$ hold along
$\partial M$.
This shows that $\omega'\in H_0^2(M)$.
On the whole, $H_0^2(M)=V^\perp$ and the orthogonal splitting is proved.\hfill$\square$\\\\
It remains to notice that, for any $\omega\in\Omega^p(M)$ with
$\omega_{|_{\partial M}}=0$, we can split
$\left(\cdot\,,\cdot\right)_V$-orthogonally $\omega=\omega_V+\hat{\omega}$,
where $\omega_V\in V$ and $\hat{\omega}\in H_0^2(M)$.
Then
\[\|\Delta\omega\|_{L^2(M)}^2=\|\omega\|_V^2=\|\omega_V\|_V^2+\|\hat{\omega}
\|_V^2\]
on the one hand, and
\[\|\nu\lrcorner d\omega\|_{L^2(\partial
M)}^2+\|\iota^*\delta\omega\|_{L^2(\partial
M)}^2=\|\nu\lrcorner
d\omega_V\|_{L^2(\partial M)}^2+\|\iota^*\delta\omega_V\|_{L^2(\partial
M)}^2\]
on the other hand, so that
\begin{eqnarray*}
\|\Delta\omega\|_{L^2(M)}^2&=&\|\Delta\omega_V\|_{L^2(M)}
^2+\|\Delta\hat {
\omega}\|_{L^2(M)}^2\\
&\geq&\|\Delta\omega_V\|_{L^2(M)}^2\\
&\geq&
q_{1,p}\cdot\left(\|\nu\lrcorner d\omega_V\|_{L^2(\partial
M)}^2+\|\iota^*\delta\omega_V\|_{L^2(\partial
M)}^2\right)\\
&=&q_{1,p}\cdot\left(\|\nu\lrcorner
d\omega\|_{L^2(\partial M)}^2+\|\iota^*\delta\omega\|_{L^2(\partial
M)}^2\right),
\end{eqnarray*}
which proves \eqref{eq:char1}.
Furthermore, the r.h.s. of
\eqref{eq:char1} is actually a minimum attained exactly by those eigenforms of
the biharmonic Steklov problem that are associated to the smallest
positive eigenvalue
$q_{1,p}$.\\

We now prove the following lemma:
\begin{lemma}
	The infimum in \eqref{eq:char2}, that we denote by $q_{1,p}'$, is a positive minimum.
\end{lemma}
{\it Proof.}  \hspace{1mm} To prove this result, we apply
the same argument as in \cite[p. 318]{FGW05}.
Namely the standard Rellich-Kondrachov compactness theorem ensures the
natural inclusion map $H^{\frac{1}{2}}(\partial M)\to L^2(\partial M)$ to be
compact.
Hence its transpose map $L^2(\partial M)\to
H^{-\frac{1}{2}}(\partial M)$ is compact, as a straightforward consequence. 
Moreover, the harmonic extension from $\partial M$ to $M$ defines a bounded
linear map
$H^{-\frac{1}{2}}(\partial M)\to L^2(M)$, see e.g. \cite[Thm. 6.6 chap.
2]{LionsMagenes68}.
Therefore the composition $L^2(\partial M)\to H^{-\frac{1}{2}}(\partial M)\to
L^2(M)$ of both maps defines a compact linear map
$E\colon L^2(\partial M)\to L^2(M)$, which already shows that $q_{1,p}'$ to be positive using only the boundedness of
the map.
Furthermore, because the image by $E$ of the unit sphere
$\mathbb{S}:=\left\{\omega\in
L^2(\partial M)\,,\,\|\omega\|_{L^2(\partial
M)}=1\right\}$ of $L^2(\partial M)$
is relatively compact in
$L^2(M)$, there exists $\tilde{\omega}\in\overline{E(\mathbb{S})}$ such that
$\|\tilde{\omega}\|_{L^2(M)}=\sup\left\{\|E\omega\|_{L^2(M)}\,|\,
\omega\in\mathbb{S} \right\}$.
By definition of the closure, there exists a sequence $(\omega_m)_m$ in
$\mathbb{S}$ such that
$E\omega_m\buil{\longrightarrow}{m\to\infty}\tilde{\omega}$ in $L^2(M)$.
But then
$\Delta(E\omega_m)\buil{\longrightarrow}{m\to\infty}\Delta\tilde{\omega}$ in
$H^{-2}(M)=(H_0^2(M))'$, so that necessarily $\Delta\tilde{\omega}=0$ holds in
$H^{-2}(M)$.
Note that, because $\Delta\tilde{\omega}=0$ and $\Delta(E\omega_m)=0$ for all
$m$, we can also claim that
$\Delta(E\omega_m)\buil{\longrightarrow}{m\to\infty}\Delta\tilde{\omega}$ in
$L^2(M)$.
By G\aa{}rding's inequality and
since both $E\omega_m\buil{\longrightarrow}{m\to\infty}\tilde{\omega}$
and $\Delta(E\omega_m)\buil{\longrightarrow}{m\to\infty}\Delta\tilde{\omega}$
in $L^2(M)$, we have that
$E\omega_m\buil{\longrightarrow}{m\to\infty}\tilde{\omega}$ in $H^2(M)$.
As a consequence, because of $\nu\lrcorner\,E\omega_m=0$ along $\partial M$ and
$\nu\lrcorner\,E\omega_m\buil{\longrightarrow}{m\to\infty
} \nu\lrcorner\,\tilde{\omega}$ in
$L^2(\partial M)$, necessarily $\nu\lrcorner\,\tilde{\omega}=0$ holds along
$\partial
M$.
Now again $E\omega_m\buil{\longrightarrow}{m\to\infty}\tilde{\omega}$
in $H^2(M)$ also implies
$\omega_m=(E\omega_m)_{|_{\partial
M}}\buil{\longrightarrow}{m\to\infty}\tilde{\omega}_{|_{\partial M}}$ in
$L^2(\partial M)$ (actually also in $H^1(\partial M)$), so that
$\|\tilde{\omega}_{|_{\partial
M}}\|_{L^2(\partial
M)}=1$.
This shows that $\tilde{\omega}\in E(\mathbb{S})$ and
hence $\|\tilde{\omega}\|_{L^2(M)}=\max\left\{\|E\omega\|_{L^2(M)}\,|\,
\omega\in\mathbb{S} \right\}$ satisfies
$\|\tilde{\omega}\|_{L^2(M)}=\frac{1}{q_{1,p}'}$.
This shows the positive r.h.s. $q_{1,p}'$ of \eqref{eq:char2} to be a minimum.\hfill$\square$

In order to finish the proof of Theorem \ref{vari_BS}, we want to show that $q_{1,p}=q_{1,p}'$.
Pick any eigenform $\omega$
associated to $q_{1,p}$.
Up to rescaling $\omega$, we may assume that $\|\nu\lrcorner
d\omega\|_{L^2(\partial M)}^2+\|\iota^*\delta\omega\|_{L^2(\partial M)}^2=1$.
Let $\hat{\omega}\in\Omega^p(M)$ be the unique solution to
$\Delta\hat{\omega}=0$ on $M$ with $\iota^*\hat{\omega}_{|_{\partial
M}}=-q_{1,p}\nu\lrcorner d\omega$ a well as
$\nu\lrcorner\hat{\omega}=q_{1,p}\iota^*\delta\omega$ along $\partial M$.
By \eqref{eq:partialintDelta},
\begin{eqnarray*}
0&=&\int_M\langle\Delta\hat{\omega},\omega\rangle\,
d\mu_g\\
&=&\int_M\langle\hat {
\omega},\Delta\omega\rangle\,d\mu_g+q_{1,p}\int_{\partial M}(|\nu\lrcorner
d\omega|^2+|\iota^*\delta\omega|^2)\,d\mu_g\\
&=&\int_M\langle\hat{
\omega},\Delta\omega\rangle\,d\mu_g+q_{1,p},
\end{eqnarray*}
so that Cauchy-Schwarz
inequality leads to
$q_{1,p}\leq\|\Delta\omega\|_{L^2(M)}\cdot\|\hat{\omega}\|_{L^2(M)}$.
Using
$q_{1,p}=\|\Delta\omega\|_{L^2(M)}^2$, we obtain
$q_{1,p}\leq\|\hat{\omega}\|_{L^2(M)}^2$.
Therefore,
\[ \frac{\|\hat{\omega}\|_{L^2(\partial
M)}^2}{\|\hat{\omega}\|_{L^2(M)}^2}=\frac{q_{1,p}^2}{\|\hat{\omega}\|_{L^2(M)}^2
} \leq q_{1,p},\]
from which $q_{1,p}'\leq q_{1,p}$ follows.
Conversely, if $\omega\in\Omega^p(M)\setminus\{0\}$ satisfies $\Delta\omega=0$
on $M$ and $\|\omega\|_{L^2(\partial
M)}^2=q_{1,p}'\|\omega\|_{L^2(M)}^2$, then let $\hat{\omega}$ be the solution
to $\Delta\hat{\omega}=\omega$ on $M$ with $\hat{\omega}_{|_{\partial
M}}=0$.
Then, again by \eqref{eq:partialintDelta}, $\omega_{|_{\partial
M}}=\iota^*\omega+\nu^\flat\wedge(\nu\lrcorner\,\omega)$ and Cauchy-Schwarz
inequality, we have
\begin{eqnarray*}
\left\|\omega\right\|_{L^2(M)}^2&=&\left(\Delta\hat{\omega},
\omega\right)_{
L^2(M) }\\
&=&\left(\hat { \omega } ,
\underbrace{\Delta\omega}_{0}\right)_{L^2(M)}+\left(\nu\lrcorner
d\hat{\omega},\iota^*\omega\right)_{L^2(\partial
M)}-\left(\iota^*\delta\hat{\omega},\nu\lrcorner\omega\right)_{L^2(\partial
M)}\\
&=&\left(\nu\lrcorner d\hat{\omega},\omega\right)_{L^2(\partial
M)}-\left(\nu^\flat\wedge\iota^*\delta\hat{\omega},\omega\right)_{L^2(\partial
M)}\\
&=&\left(\nu\lrcorner\,d\hat{\omega}-\nu^\flat\wedge\iota^*\delta\hat{\omega},
\omega\right)_{L^2(\partial
M)}\\
&\leq&\left\|\nu\lrcorner\,d\hat{\omega}-\nu^\flat\wedge\iota^*\delta\hat{\omega
}\right\|_{L^2(\partial M)}\cdot\left\|\omega\right\|_{L^2(\partial M)}.
\end{eqnarray*}
But since
$\left(\nu\lrcorner\,d\hat{\omega},\nu^\flat\wedge\iota^*\delta\hat{\omega
}\right)_{L^2(\partial M)}=0$, we have
\[\left\|\nu\lrcorner\,d\hat{\omega}-\nu^\flat\wedge\iota^*\delta\hat{\omega
}\right\|_{L^2(\partial
M)}^2=\left\|\nu\lrcorner\,d\hat{\omega}\right\|_{L^2(\partial
M)}^2+\left\|\nu^\flat\wedge\iota^*\delta\hat{\omega
}\right\|_{L^2(\partial
M)}^2=\left\|\nu\lrcorner\,d\hat{\omega}\right\|_{L^2(\partial
M)}^2+\left\|\iota^*\delta\hat{\omega
}\right\|_{L^2(\partial
M)}^2,\]
so that
$\left\|\omega\right\|_{L^2(M)}^2
\leq \left(\left\|\nu\lrcorner
d\hat{\omega}\right\|_{L^2(\partial
M)}^2+\left\|\iota^*\delta\hat{\omega}\right\|_{L^2(\partial
M)}^2\right)^{\frac{1}{2}}\cdot\left\|\omega\right\|_{L^2(\partial M)}$.
Therefore,
\begin{eqnarray*}
\frac{\|\Delta\hat{\omega}\|_{L^2(M)}^2}{\|\nu\lrcorner
d\hat{\omega}\|_{L^2(\partial
M)}^2+\|\iota^*\delta\hat{\omega}\|_{L^2(\partial
M)}^2}&=&\frac{\|\omega\|_{L^2(M)}^2}{\|\nu\lrcorner
d\hat{\omega}\|_{L^2(\partial
M)}^2+\|\iota^*\delta\hat{\omega}\|_{L^2(\partial
M)}^2}\\
&\leq&\frac{\|\omega\|_{L^2(M)}^2\cdot\left\|\omega\right\|_{L^2(\partial
M)}^2}{ \|\omega\|_{ L^2(M)
}^4}\\
&=&\frac{\left\|\omega\right\|_{L^2(\partial
M)}^2}{ \|\omega\|_{ L^2(M) }^2}\\
&=&q_{1,p}',
\end{eqnarray*}
from
which $q_{1,p}'\geq q_{1,p}$
follows.
On the whole, we deduce that
$q_{1,p}=q_{1,p}'$, as we claimed.
Moreover, the $p$-form $\hat{\omega}$ defined above, because it now minimizes $q_{1,p}$, must be an eigenform of the biharmonic Steklov operator associated to the eigenvalue $q_{1,p}$.
As a consequence, $\omega=\Delta\hat{\omega}$ where $ \hat{\omega}$ is an eigenform of the biharmonic Steklov operator associated to the eigenvalue $q_{1,p}$.
\hfill$\square$

\section{Eigenvalues of the biharmonic Steklov operator}\label{s:ebs}
In this section, we will establish some eigenvalue estimates for the first eigenvalue of the biharmonic Steklov operator defined in the previous section.

As before, we will consider a compact Riemannian manifold $(M^n,g)$ with smooth boundary $\partial M$. Notice first the following fact:

\begin{lemma} The biharmonic Steklov operator is preserved by the Hodge star operator $*_M$ on $M$.
\end{lemma}
{\it Proof.} 
We only need to check that the last two boundary conditions in \eqref{biharmonic-Steklov} are preserved. 
For this purpose, using the equality $\iota^*(*_M\alpha)=*_{\partial M}(\nu\lrcorner\alpha)$ for any form $\alpha$, we compute, for any solution $\omega$ of degree $p$ to problem \eqref{biharmonic-Steklov},
\begin{eqnarray*}
\nu\lrcorner
\Delta(*_M\omega)&=&\nu\lrcorner(*_M\Delta\omega)\\
&=&(-1)^p*_{\partial
M}\iota^*(\Delta\omega)\\
&=&(-1)^p q*_{\partial M}(\nu\lrcorner d\omega)\\
&=&(-1)^p
q\iota^*(*_Md\omega)\\
&=&-q\iota^*(\delta*_M\omega).
\end{eqnarray*}
In the last equality, we use the fact that $*_Md=(-1)^{p-1}\delta*_M$ on $p$-forms. For the other boundary condition, we have
$$\iota^*(\Delta *_M\omega)=\iota^*(*_M\Delta\omega)=*_{\partial M}(\nu\lrcorner\Delta\omega)=-q*_{\partial M}(\iota^*\delta\omega)=(-1)^pq\nu\lrcorner*_M\delta\omega=q\nu\lrcorner d(*_M\omega).$$
Also here we use the fact that $d*_M=(-1)^{p}*_M\delta$ on $p$-forms. This finishes the proof.
\hfill$\square$
\begin{remark}
\rm As a direct consequence of the invariance of the biharmonic Steklov operator by the Hodge star operator is that $q_{i,p}=q_{i,n-p}$ for any $i\geq 1$ and $p\leq n$.
\end{remark}

In the following, we recall the estimate stated by S. Raulot and A. Savo in \cite{RaulotSavo2015} for subharmonic functions that we will use in order to get a lower bound of the first eigenvalue $q_{1,p}$. 
Let $(M^n,g)$ be a compact Riemannian manifold with smooth boundary such that the Ricci curvature of $M$ satisfies ${\rm Ric}_M\geq (n-1)K$ and the mean curvature of the boundary satisfies $H\geq H_0$ for some real numbers $K$ and $H_0$. 
Let $R$ be the inner radius of the manifold $M$, that is
$$R={\rm max}\{{\rm dist}(x,\partial M)|\, x\in M\},$$
and $\Theta(r)=(s'_K(r)-H_0s_K(r))^{n-1}$  for all $r$, where the function $s_K$ is being
given by
\[
s_K(r):=
\left\{
\begin{array}{lll}\medskip
\frac{1}{\sqrt{K}}\,\sin(r\sqrt{K}) && \text{if}\,\, K>0, \\
\medskip
r && \text{if}\,\, K=0,\\
\medskip
\frac{1}{\sqrt{|K|}}\,\sinh(r\sqrt{|K|}) && \text{if}\,\, K<0.
\end{array}
\right.
\]
It was shown in \cite[Prop. 14]{RaulotSavo2015} (see also \cite[Thm. A]{Kasue}) that the function $\Theta$ is smooth and positive on $[0,R[$ and $\Theta(R)=0$ when $M$ is a geodesic ball in $M_K$, the space form of sectional curvature $K$. 
The following result was proved in  \cite{RaulotSavo2015}:

\begin{thm} \cite[Thm. 10]{RaulotSavo2015} Let $(M^n,g)$ be a compact Riemannian manifold with smooth boundary. Assume that the Ricci curvature of $M$ satisfies ${\rm Ric}_M\geq (n-1)K$ and the mean curvature $H\geq H_0$ for some real numbers $K$ and $H_0$. If $h$ is a non-trivial, nonnegative subharmonic function on $M$ (i.e. $\Delta h\leq 0$ on $M$), then
$$\frac{\int_{\partial M} hd\mu_g}{\int_M h d\mu_g}\geq \frac{1}{\int_0^R \Theta(r) dr}.$$
\end{thm}

Using this result and the Bochner formula $\Delta=\nabla^*\nabla+\mathfrak{B}^{[p]}$ on $p$-forms, we  prove the following:

\begin{thm} Let $(M^n,g)$ be a compact Riemannian manifold with smooth boundary. Assume that the Ricci curvature of $M$ satisfies ${\rm Ric}_M\geq (n-1)K$ and the mean curvature $H\geq H_0$ for some real numbers $K$ and $H_0$. Assume also that the Bochner operator $\mathfrak{B}^{[p]}$ is nonnegative for some $p$. Then, the inequality
\begin{equation}\label{eq:minq1Theta}
q_{1,p}\geq\frac{1}{\int_0^R\Theta(r)\,dr}
\end{equation}
holds.
\end{thm}

{\it Proof.} Applying the Bochner formula to $\Delta\omega$ where $\omega$ is a $p$-eigenform of the biharmonic Steklov operator associated with $q_{1,p}$, we get after taking the pointwise scalar product with $\Delta\omega$ that
$$0=\langle\Delta^2\omega,\Delta\omega\rangle=|\nabla\Delta\omega|^2+\frac{1}{2}\Delta(|\Delta\omega|^2)+\langle \mathfrak{B}^{[p]}\Delta\omega,\Delta\omega\rangle.$$
Since $\mathfrak{B}^{[p]}$ is nonnegative, we deduce that $\Delta(|\Delta\omega|^2)$ is nonpositive or equivalently the function $h:=|\Delta\omega|^2$ is subharmonic. Therefore, by the previous theorem, we can say that
$$\frac{\int_{\partial M} |\Delta\omega|^2 d\mu_g}{\int_M |\Delta\omega|^2 d\mu_g}\geq \frac{1}{\int_0^R \Theta(r) dr}.$$
Now, Characterization \eqref{eq:char2} gives the result and finishes the proof of the theorem.
\hfill$\square$

\begin{remark}
\rm Depending on the sign of $K$ and $H_0$, we can estimate explicitly $\int_0^R\Theta(r) dr$ in terms of $R$ and $H_0$, as in \cite[Thms. 12 \& 13]{RaulotSavo2015}. Therefore, one can deduce several estimates for $q_{1,p}$ in terms of $R$ and $H_0$.
\end{remark}


We will now provide an estimate for the first eigenvalue of problem \eqref{biharmonic-Steklov} on manifolds carrying parallel forms and study the limiting case of the estimate.
Recall that a harmonic domain is a compact Riemannian manifold $(M^n,g)$ with smooth boundary $\partial M$ supporting a solution to the Serrin problem \eqref{eq:scalarSerrin}.
We have the following result:
\begin{thm}\label{thm:harmonicdomain}
Let $(M^n,g)$ be a compact Riemannian manifold with smooth boundary. Assume that $M$ supports a non-trivial parallel $p$-form $\omega_0$ for some $p=0,\ldots,n$.
Then
\begin{equation} \label{estharmdom}
q_{1,p} \le \frac{{\rm Vol}(\partial M)}{{\rm Vol}(M)}.\end{equation}
Moreover, if equality holds in \eqref{estharmdom}, then $f\cdot\omega_0$ is an eigenform associated to $q_{1,p}$, where $f$ is the solution of \eqref{eq:scalarSerrin} and therefore $M$ must be a harmonic domain (and hence a Euclidean ball if $M\subset\R^n$).\\
Conversely, if $M$ is a harmonic domain, then $\frac{{\rm Vol}(\partial M)}{{\rm Vol}(M)}$ is an eigenvalue of the biharmonic Steklov problem \eqref{biharmonic-Steklov}.
\end{thm}
{\it Proof.} As $\omega_0$ is a parallel form, we can assume that $|\omega_0|=1$. By using the variational characterization \eqref{eq:char2}, we obtain that
$$q_{1,p} \le \dfrac{\|\omega_0\|_{L^2(\partial M)}^2}{\|\omega_0\|_{L^2(
M)}^2}=\dfrac{{\rm Vol}(\partial M)}{{\rm Vol}(M)}.$$
If equality occurs in \eqref{estharmdom}, then $\omega_0=\Delta \omega$ for some eigenform $\omega$ associated with $q_{1,p}$ by Theorem \ref{vari_BS}.
Now Proposition \ref{prop:serrinformes} implies that $M$ carries a solution $f$ to the Serrin problem \eqref{eq:scalarSerrin} and that $\omega=f\cdot\omega_0$.
To check the converse, we take a function $f$ solution to the Serrin problem \eqref{eq:scalarSerrin}, then using again Proposition \ref{prop:serrinformes} the $p$-form $\omega:=f\cdot \omega_0$ is an eigenform of \eqref{biharmonic-Steklov} associated with the eigenvalue $q=\frac{1}{c}$, where $c=\dfrac{{\rm Vol}( M)}{{\rm Vol}(\partial M)}$.
\hfill$\square$

\begin{remark}\label{rem:harmonicformofconstantlength}
\rm
If a compact Riemannian manifold $M$ with smooth boundary carries a nontrivial \emph{harmonic form of constant length} $\omega_0$, then \eqref{estharmdom} remains valid.
Moreover, if \eqref{estharmdom} is an equality, then there still exists an eigenform $\omega$ to the biharmonic Steklov operator on $p$-forms such that $\Delta\omega=\omega_0$, nevertheless it is no more true in general that $M$ must be a harmonic domain and that $\omega=f\cdot\omega_0$ for some solution $f$ to the scalar Serrin problem.
\end{remark}

Next we compare the first
eigenvalues of the biharmonic Steklov operator for successive degrees, when the
manifold $M$ is a domain in $\R^n$ or $\mathbb{S}^n$.
We first notice that, if $f$ is any eigenfunction to the scalar biharmonic Steklov problem \eqref{biharmonic-Steklovbis}, then for any parallel $p$-form $\omega_0$ on $M$ the form $f\cdot\omega_0$ is an eigenform to the biharmonic Steklov problems on $p$-forms and associated to the \emph{same eigenvalue}.
Therefore, for every eigenvalue $q$ of the scalar biharmonic Steklov operator, we have an embedding
\[ \ker(\mathrm{BS}_0-q)\otimes\mathrm{P}_p\hookrightarrow\ker(\mathrm{BS}_p-q),
\]
where $\ker(\mathrm{BS}_j-q)$ denotes the eigenspace for the biharmonic Steklov operator on $j$-forms and associated to the eigenvalue $q$, and $\mathrm{P}_p$ denotes the space of parallel $p$-forms on $M$.
When $M\subset\R^n$, then conversely for any $\omega\in\ker(\mathrm{BS}_p-q)$, there exists a parallel $p$-form $\omega_0$ on $M$ with $|\omega_0|=1$ and $\langle\omega,\omega_0\rangle\neq0$ (non identically vanishing) on $M$.
But then $f$ can be easily shown to lie in $\ker(\mathrm{BS}_0-q)$.
This shows that, when $M\subset \R^n$, both $0$- and $p$-biharmonic Steklov eigenvalues coincide, their multiplicities being ignored.

In what follows, we assume that $(M^n,g)$ is isometrically immersed into the Euclidean
space $\R^{n+m}$. For any given smooth  normal vector field $N$ to $M$, we
denote by $I\!I_N$ the associated Weingarten map, that is, the
endomorphism field of $TM$ defined by
\[\langle I\!I_N(X),Y\rangle=\langle N,I\!I(X,Y)\rangle\]
for all $X,Y$ tangent to $M$, where $I\!I$ is the second fundamental form of the
immersion. Recall that any endomorphism $A$ of $TM$ can be extended to the set of differential $p$-forms on $M$ as follows: For any $p$-form $\omega$ on $M$, we define
\begin{equation}\label{extensionweingarten}
A^{[p]}\omega(X_1,\cdots, X_p)=\sum_{i=1}^p \omega(X_1,\cdots, A(X_i),\cdots, X_p),
\end{equation}
for all $X_1,\cdots, X_p$ vector fields in $TM$. In particular, this applies to $I\!I_N$ for all $N\in T^\perp M$.  The following lemma is technical but will be useful for the comparison.

\begin{lemma} \label{lem:inequalitysteklov} Let $(M^n,g)$ be a compact Riemannian manifold with smooth boundary $\partial M$. Assume that $M$ is isometrically immersed into the Euclidean space $\R^{n+m}$. Let $\omega$ be any $p$-eigenform of the biharmonic Steklov operator. Then we have
\begin{eqnarray}\label{eq:inequalityimmersion}
p q_{1,p-1} \int_{\partial M}(|\nu\lrcorner d\omega|^2+|\iota^*(\delta\omega)|^2)d\mu_g\leq
 \sum_{i=1}^{n+m}\int_M|-2I\!I_{n\widetilde{H}}(\partial_{x_i}^T)\lrcorner\omega+2\sum_{a=1}^m I\!I_{f_a}(\partial_{x_i}^T) \lrcorner I\!I_{f_a}^{[p]}\omega\nonumber\\-2\sum_{s=1}^n (I\!I_{\partial_{x_i}^\perp} e_s) \lrcorner (   \nabla_{e_s} \omega )-n\left(d(\langle \widetilde{H},\partial_{x_i}^\perp\rangle)\right)\lrcorner \omega+\partial_{x_i}^T\lrcorner \Delta \omega|^2 d\mu_g.\nonumber\\
\end{eqnarray}
Here $\{e_1,\cdots,e_n\}$ and $\{f_1,\cdots,f_m\}$ are respectively orthonormal bases of $TM$ and $T^\perp M$ and $\widetilde{H}$ is the mean curvature field of the immersion.
\end{lemma}

{\it Proof.} Let $\omega$ be any eigenform of the biharmonic Steklov problem
associated with $q_{1,p}$. For each $i=1,\cdots, n+m$, the unit parallel vector field $\partial_{x_i}$ on
$\R^{n+m}$ splits into
$\partial_{x_i}=(\partial_{x_i})^T+(\partial_{x_i})^\perp$ where
$(\partial_{x_i})^T$ is the tangent part in $TM$ and $(\partial_{x_i})^\perp$ is
the orthogonal one in $T^\perp M$.
We consider the $(p-1)$-form $(\partial_{x_i})^T\lrcorner \omega$ on $M$ which
clearly vanishes on $\partial M$. By
applying to it the variational characterization \eqref{eq:char1},
we get, for each $i$,
\begin{equation}\label{eq:testpartialxi1}
q_{1,p-1} \int_{\partial M}\left(|\nu \lrcorner d((\partial_{x_i})^T\lrcorner
\omega)|^2 +|\iota^*  \delta ( (\partial_{x_i})^T\lrcorner
\omega)|^2  \right)d\mu_g\le
\int_M |\Delta((\partial_{x_i})^T\lrcorner
\omega)|^2d\mu_g.\end{equation}
Now we want to sum over $i=1,\cdots,n+m$. We first begin with the l.h.s. 
Recall the Cartan
formula: $\mathcal L_X \omega= d(X\lrcorner \omega)+X\lrcorner d\omega$, for any vector field $X$ on $M$. Using this formula, we have for each $i$,
\begin{eqnarray}\label{eq:splitLie}
d(\partial_{x_i}^T\lrcorner\omega)&=&\mathcal L_{\partial_{x_i}^T}\omega-\partial_{x_i}
^T\lrcorner d\omega\nonumber\\
&=&\nabla_{\partial_{x_i}^T}\omega+I\!I_{\partial_{x_i}^\perp}^{[p]}\omega-\partial_{x_i}
^T\lrcorner d\omega.\end{eqnarray}
In the last equality, we used the splitting of the Lie derivative in terms of the connection as follows: $\mathcal L_{X^T}\omega=\nabla_{X^T}\omega+I\!I_{X^\perp}^{[p]}\omega$, for a parallel vector field $X \in \R^{n+m}$, proved in \cite[Eq. (4.3) p. 337]{GS}. Since $\omega=0$ on $\partial M$, we have that for any $X\in T\partial M$ \cite[Eq. (23)]{RaulotSavo2011}
$$\nu\lrcorner \nabla_X\omega=\nabla_X^{\partial M}(\nu\lrcorner\omega)+S(X)\lrcorner(\iota^*\omega)=0.$$
Here $S$ denotes the second fundamental form of the boundary. Therefore, we deduce that
\begin{eqnarray}\label{eq:nudx}
\nu \lrcorner d((\partial_{x_i})^T\lrcorner
\omega)&=&\nu \lrcorner \nabla_{\partial_{
x_i}^T}\omega +\partial_{x_i}
^T\lrcorner (\nu \lrcorner d\omega)\nonumber\\
&=&g((\partial_{x_i})^T,\nu)\nu\lrcorner\nabla_\nu\omega +\partial_{x_i}
^T\lrcorner (\nu \lrcorner d\omega)\nonumber\\
&=&-g(\partial_{x_i},\nu)\iota^*(\delta\omega)+\partial_{x_i}
^T\lrcorner (\nu \lrcorner d\omega).
\end{eqnarray}
In the last equality, we use the identity \cite[Lem. 18]{RaulotSavo2011}
$$\nu\lrcorner\nabla_\nu\omega=\delta^{\partial M}(\iota^*\omega)-\iota^*(\delta\omega)-S^{[p-1]}(\nu\lrcorner\omega)+(n-1)H\nu\lrcorner\omega=-\iota^*(\delta\omega).$$
Independently, by a straightforward computation, we check that, for any $p$-form $\alpha$,
$$
\sum_{i=1}^{n+m}\partial_{x_i}^T\wedge (\partial_{x_i}^T\lrcorner \alpha)=\sum_{i=1}^n e_i \wedge (e_i \lrcorner \alpha)=p\alpha.
$$
As a consequence, we obtain for any $p$-forms $\alpha$ and $\beta$ on $M$
\begin{equation}\label{eq:scalinnerpro}
\sum_{i=1}^{n+m}\langle \partial_{x_i}^T\lrcorner \alpha, \partial_{x_i}^T\lrcorner \beta \rangle=p\langle \alpha, \beta \rangle.
\end{equation}
We take the norm of \eqref{eq:nudx} and sum over $i$ to get
\begin{eqnarray}\label{eq:nud}
 \sum_{i=1}^{n+m}|\nu \lrcorner d((\partial_{x_i})^T\lrcorner
\omega)|^2 &\bui{=}{\eqref{eq:scalinnerpro}}& |\iota^*(\delta\omega)|^2 +p|\nu \lrcorner d\omega|^2 - 2 \sum_{i=1}^{n+m} g(\partial_{x_i},\nu)\langle \iota^*(\delta\omega) ,  \partial_{x_i}
^T\lrcorner ( \nu \lrcorner d\omega)\rangle \nonumber\\
&=&| \iota^*\delta\omega|^2 +p|\nu \lrcorner d\omega|^2-2\langle \iota^*\delta\omega, \nu^T\lrcorner ( \nu \lrcorner d\omega)\rangle \nonumber\\
&=& | \iota^*\delta\omega|^2 +p|\nu \lrcorner d\omega|^2.
\end{eqnarray}
Here we notice that by the Gau\ss{} formula and the fact that $\partial_{x_i}$ is parallel in $\R^{n+m}$, we have $\nabla_X \partial_{x_i}^T=I\!I_{\partial_{x_i}^\perp}(X)$, so that
\begin{eqnarray*}
\delta( \partial_{x_i}^T \lrcorner
\omega)&=&-\sum_{j=1}^ne_j\lrcorner\,\nabla_{e_j}(\partial_{x_i}^T \lrcorner
\omega)\\
&=&-\sum_{j=1}^ne_j\lrcorner\,(\nabla_{e_j}\partial_{x_i}^T\lrcorner\,\omega+\partial_{x_i}^T\lrcorner\,\nabla_{e_j}\omega)\\
&=&-\sum_{j=1}^ne_j\lrcorner\,I\!I_{\partial_{x_i}^\perp}(e_j)\lrcorner\,\omega+\sum_{j=1}^n\partial_{x_i}^T\lrcorner\,(e_j\lrcorner\,\nabla_{e_j}\omega)\\
&=&-\sum_{j,k=1}^n\langle I\!I_{\partial_{x_i}^\perp}(e_j),e_k\rangle e_j\lrcorner\, e_k\lrcorner\,\omega-\partial_{x_i}^T\lrcorner\,\delta\omega\\
&=&-\partial_{x_i}^T\lrcorner\,\delta\omega
\end{eqnarray*}
because of the expression $\langle I\!I_{\partial_{x_i}^\perp}(e_j),e_k\rangle =\langle I\!I(e_j,e_k),\partial_{x_i}^\perp\rangle$ being symmetric in $j,k$.
Thus $\delta( \partial_{x_i}^T \lrcorner
\omega)=- \partial_{x_i}^T \lrcorner  \delta\omega$ for any $i$.
Using also the fact that $\nu\lrcorner\,\delta(\partial_{x_i}^T \lrcorner
\omega)=-\delta^{\partial M}(\nu\lrcorner\,\partial_{x_i}^T \lrcorner
\omega)=0$, we get that
\begin{equation}\label{eq:idelta}
\sum_{i=1}^{n+m} | \iota^*\delta( \partial_{x_i}^T \lrcorner
\omega)|^2=\sum_{i=1}^{n+m} | \delta( \partial_{x_i}^T \lrcorner
\omega)|^2 \bui{=}{\eqref{eq:scalinnerpro}} (p-1)| \iota^*\delta\omega|^2.
\end{equation}
Hence adding Equations \eqref{eq:nud} and \eqref{eq:idelta} allows to find the l.h.s. of Inequality \eqref{eq:inequalityimmersion}.
We are now going to estimate the term
$\displaystyle\sum_{i=1}^{n+m}
 |\Delta(\partial_{x_i}^T\lrcorner
\omega)|^2$ in \eqref{eq:testpartialxi1}. Taking the divergence of \eqref{eq:splitLie} and the differential of the identity $\delta( \partial_{x_i}^T \lrcorner
\omega)=- \partial_{x_i}^T \lrcorner  \delta\omega$ along with the Cartan formula and the decomposition of the Lie derivative as in \eqref{eq:splitLie}, we get that
\begin{equation}\label{Delta-test1}
\Delta (\partial_{x_i}^T\lrcorner\omega)=[\delta,\nabla_{\partial_{
x_i}^T}](\omega)+\delta(I\!I_{\partial_{x_i}^\perp}^{[p]}\omega)-I\!I^{[p-1]}_{\partial_{x_i}^\perp}(\delta\omega)+\partial_{x_i}
^T\lrcorner \Delta \omega.
\end{equation}

In the following, we will compute each term of \eqref{Delta-test1} separately. 
First, take an orthonormal frame $\{e_1,\cdots,e_n\}$ of $TM$ such that $\nabla e_i=0$ at some point.
Then we have, for any vector field $X$ on $M$,
\begin{eqnarray}\label{firstbracket}
[\delta,\nabla_{X}](\omega) &=& \delta(\nabla_X \omega) - \nabla_X \delta \omega \nonumber \\
&=& -\sum_{s=1}^n e_s \lrcorner \nabla_{e_s} \nabla_X \omega - \nabla_X \delta \omega \nonumber \\
&=&-\sum_{s=1}^n e_s \lrcorner \left(R(e_s,X) \omega +\nabla_X \nabla_{e_s}  \omega + \nabla_{[e_s,X]}\omega\right) - \nabla_X \delta \omega \nonumber\\
&=&-\sum_{s=1}^n e_s \lrcorner R(e_s,X) \omega +\nabla_X\delta\omega -\sum_{s=1}^n e_s \lrcorner \nabla_{\nabla_{e_s}X}\omega - \nabla_X \delta \omega \nonumber \\
&=& -\sum_{s=1}^n e_s \lrcorner R(e_s,X) \omega -\sum_{s=1}^n\ e_s \lrcorner \nabla_{\nabla_{e_s}X}\omega.
\end{eqnarray}
Now we use the fact that for any tensor field $A$, any vector field $X$ and any $p$-form $\alpha$ on $M$,  $\nabla_X A^{[p]}=(\nabla_X A)^{[p]}$ and $A^{[p-1]}(X\lrcorner \alpha)=X\lrcorner A^{[p]}\alpha-A(X)\lrcorner\alpha$, which both can be proved by a straightforward computation.
Thus for any $N \in T^\perp M$, we can write
\begin{eqnarray} \label{eq:secondbracket}
\delta (I\!I^{[p]}_N \omega) -I\!I^{[p-1]}_N (\delta \omega) \nonumber  &=& -\sum_{s=1}^n e_s \lrcorner \nabla_{e_s} (I\!I^{[p]}_N \omega) -I\!I^{[p-1]}_N (\delta \omega) \nonumber \\
&=& -\sum_{s=1}^n e_s \lrcorner ( \nabla_{e_s} I\!I_N)^{[p]} (\omega)-\sum_{s=1}^n e_s \lrcorner  I\!I_N^{[p]}(  \nabla_{e_s} \omega ) -I\!I^{[p-1]}_N (\delta \omega)   \nonumber \\
&=& -\sum_{s=1}^n (\nabla_{e_s} I\!I_N)^{[p-1]} ( e_s \lrcorner \omega)-  \sum_{s=1}^n (( \nabla_{e_s} I\!I_N)(e_s))\lrcorner \omega-\sum_{s=1}^n I\!I^{[p-1]}_N (e_s \lrcorner \nabla_{e_s} \omega)\nonumber\\&&-\sum_{s=1}^n (I\!I_N e_s) \lrcorner (   \nabla_{e_s} \omega ) -I\!I^{[p-1]}_N (\delta \omega) \nonumber \\
&=& -\sum_{s=1}^n (\nabla_{e_s} I\!I_N)^{[p-1]} ( e_s \lrcorner \omega)+\delta(I\!I_N)\lrcorner \omega-\sum_{s=1}^n (I\!I_N e_s) \lrcorner (   \nabla_{e_s} \omega).
\end{eqnarray}
By taking $X=\partial_{x_i}^T$ in \eqref{firstbracket} and $N=\partial_{x_i}^\perp$ in \eqref{eq:secondbracket}, Equation \eqref{Delta-test1} reduces after using the fact that $\nabla\partial_{x_i}^{T}= I\!I_{ \partial_{x_i}^{\perp}}$ as a consequence of the parallelism of the vector field $\partial_{x_i}$, to
\begin{eqnarray}\label{Delta-test2}
\Delta (\partial_{x_i}^T\lrcorner\omega) &=&  -\sum_{s=1}^n e_s \lrcorner R(e_s,\partial_{x_i}^T) \omega -2\sum_{s=1}^n (I\!I_{\partial_{x_i}^\perp} e_s) \lrcorner (   \nabla_{e_s} \omega ) \nonumber \\
&&-\sum_{s=1}^n (\nabla_{e_s} I\!I_{\partial_{x_i}^\perp})^{[p-1]} ( e_s \lrcorner \omega)+\delta(I\!I_{\partial_{x_i}^\perp})\lrcorner \omega+\partial_{x_i}^T\lrcorner \Delta \omega.
\end{eqnarray}
We proceed in the computation of \eqref{Delta-test2}. Let us compute the term $\sum_{s=1}^n (\nabla_{e_s} I\!I_{\partial_{x_i}^\perp})^{[p-1]} ( e_s \lrcorner \omega)$. Using the fact that $A^{[p]}=\sum_{l=1}^n e_l\wedge \left(A(e_l)\lrcorner\right)$ for any symmetric tensor $A$, we compute with the help of Equation \eqref{eq:delta} in the appendix
\begin{eqnarray}\label{eq:nablaII}
\sum_{s=1}^n (\nabla_{e_s} I\!I_{\partial_{x_i}^\perp})^{[p-1]} ( e_s \lrcorner \omega)&=&\sum_{l,s=1}^n e_l\wedge \left((\nabla_{e_s} I\!I_{\partial_{x_i}^\perp})(e_l)\lrcorner e_s\lrcorner\omega\right)\nonumber \\
&=& \sum_{l,s=1}^ne_l\wedge \left((\nabla_{e_l} I\!I_{\partial_{x_i}^\perp})(e_s)\lrcorner e_s\lrcorner \omega\right)-\sum_{l,s=1}^ne_l\wedge  \left(I\!I_{(\nabla^{\R^{n+m}}_{e_l}\partial_{x_i}^\perp)^\perp})(e_s)\lrcorner e_s\lrcorner \omega\right)\nonumber\\
&&+\sum_{l,s=1}^ne_l\wedge  \left(I\!I_{(\nabla^{\R^{n+m}}_{e_s}\partial_{x_i}^\perp)^\perp})(e_l)\lrcorner e_s\lrcorner \omega\right).
\end{eqnarray}
The first two sums vanish identically, since $\sum_s A(e_s)\lrcorner e_s\lrcorner\omega=0$ for any symmetric endomorphism $A$. Hence with the help of Equation \eqref{eq:nablaxiperp}, Equality \eqref{eq:nablaII} reduces to
\begin{eqnarray}\label{eq:nablaii}
\sum_{s=1}^n (\nabla_{e_s} I\!I_{\partial_{x_i}^\perp})^{[p-1]} ( e_s \lrcorner \omega)&=&-\sum_{l,s=1}^ne_l\wedge  \left(I\!I_{I\!I(e_s,\partial_{x_i}^T)}(e_l)\lrcorner e_s\lrcorner \omega\right)\nonumber\\
&=&-\sum_{l,s=1}^n\sum_{a=1}^m \langle I\!I(e_s,\partial_{x_i}^T),f_a\rangle e_l\wedge  \left(I\!I_{f_a}(e_l)\lrcorner e_s\lrcorner \omega\right)\nonumber\\
&=&-\sum_{l=1}^n\sum_{a=1}^m  e_l\wedge  \left(I\!I_{f_a}(e_l)\lrcorner I\!I_{f_a}(\partial_{x_i}^T)\lrcorner \omega\right)=-\sum_{a=1}^m I\!I_{f_a}^{[p-1]}(I\!I_{f_a}(\partial_{x_i}^T)\lrcorner\omega)\nonumber\\
\end{eqnarray}
where $\{f_1,\cdots,f_m\}$ is an orthonormal frame of $T^\perp M$. In the last equality, we used again the expression $A^{[p]}=\sum_l e_l\wedge (A(e_l)\lrcorner)$ for any symmetric endomorphism $A$. Hence after replacing Equations \eqref{Sec-term-delta}, proved in the appendix, and \eqref{eq:nablaii} into Equation \eqref{Delta-test2}, we finally  get
\begin{eqnarray}\label{eq:Delta}
\Delta (\partial_{x_i}^T\lrcorner\omega) &=&  -\sum_{s=1}^n e_s \lrcorner R(e_s,\partial_{x_i}^T) \omega -2\sum_{s=1}^n (I\!I_{\partial_{x_i}^\perp} e_s) \lrcorner (   \nabla_{e_s} \omega ) \nonumber \\
&&+\sum_{a=1}^m
  I\!I_{f_a}^{[p-1]}(I\!I_{f_a}(\partial_{x_i}^T) \lrcorner \omega)+\left(\sum_{a=1}^m
  I\!I^2_{f_a}(\partial_{x_i}^T)-n d(\langle \widetilde{H},\partial_{x_i}^\perp\rangle)-n I\!I_{\widetilde{H}}(\partial_{x_i}^T)\right)\lrcorner \omega+\partial_{x_i}^T\lrcorner \Delta \omega  \nonumber\\
  &=&-\sum_{s=1}^n e_s \lrcorner R(e_s,\partial_{x_i}^T) \omega -2\sum_{s=1}^n (I\!I_{\partial_{x_i}^\perp} e_s) \lrcorner (   \nabla_{e_s} \omega ) \nonumber \\
&&+\sum_{a=1}^m
  I\!I_{f_a}(\partial_{x_i}^T) \lrcorner I\!I_{f_a}^{[p]}\omega+\left(-n d(\langle \widetilde{H},\partial_{x_i}^\perp\rangle)-n I\!I_{\widetilde{H}}(\partial_{x_i}^T)\right)\lrcorner \omega+\partial_{x_i}^T\lrcorner \Delta \omega  \nonumber\\
  &=&-I\!I_{n\widetilde{H}}(\partial_{x_i}^T)\lrcorner\omega+ \sum_{a=1}^mI\!I_{f_a}(\partial_{x_i}^T)\lrcorner I\!I_{f_a}^{[p]}\omega-2\sum_{s=1}^n (I\!I_{\partial_{x_i}^\perp} e_s) \lrcorner (   \nabla_{e_s} \omega )  \nonumber \\
  &&+\sum_{a=1}^m
  I\!I_{f_a}(\partial_{x_i}^T) \lrcorner I\!I_{f_a}^{[p]}\omega+\left(-n d(\langle \widetilde{H},\partial_{x_i}^\perp\rangle)-n I\!I_{\widetilde{H}}(\partial_{x_i}^T)\right)\lrcorner \omega+\partial_{x_i}^T\lrcorner \Delta \omega  \nonumber\\
  &=&-2I\!I_{n\widetilde{H}}(\partial_{x_i}^T)\lrcorner\omega+2\sum_{a=1}^mI\!I_{f_a}(\partial_{x_i}^T)\lrcorner I\!I_{f_a}^{[p]}\omega-2\sum_{s=1}^n (I\!I_{\partial_{x_i}^\perp} e_s) \lrcorner (   \nabla_{e_s} \omega )-n\big(d\big(\langle \widetilde{H},\partial_{x_i}^\perp\rangle)\big)\lrcorner \omega \nonumber\\
  &&+\partial_{x_i}^T\lrcorner \Delta \omega.
  \end{eqnarray}
In the second equality, we used again the relation $A^{[p-1]}(X\lrcorner \alpha)=X\lrcorner A^{[p]}\alpha-A(X)\lrcorner\alpha$ for any $p$-form $\alpha$ and $X\in TM$, and in the third equality, we used Proposition \ref{prop:curvature} in the appendix. Equation \eqref{eq:Delta}, along with \eqref{eq:nud} and \eqref{eq:idelta}, gives the result, by using Inequality \eqref{eq:testpartialxi1}.
\hfill$\square$

In general, it is difficult to control all the terms in Inequality \eqref{eq:inequalityimmersion} in order to compare $q_{1,p-1}$ and $q_{1,p}$. 
Therefore, we shall restrict ourselves to the case when $M$ is a domain in $\mathbb{S}^n$. 
We have the following result:

\begin{thm}\label{est_sphere} Let $(M^n,g)$ be a compact Riemannian manifold with smooth boundary.
If $M$ is a domain in $\mathbb{S}^n$, then we have that
\begin{equation}\label{eq:ineqq1psphere}
q_{1,p-1}+ (n-p) q_{1,p+1}< C_{p,n} q_{1,p},
\end{equation}
where $C_{p,n}$ is some constant that depends on $p$ and $n$ and whose explicit expression is given in \eqref{eq:c(p,n)}.
\end{thm}

{\it Proof.}
We consider the
 isometric immersion $M\subset \mathbb{S}^n\hookrightarrow \R^{n+1}$. 
In this
case, we have that $m=1$, the orthonormal basis of $T^\perp M$ reduces to
the inward unit vector field $\tilde\nu=-\sum_{i=1}^{n+1} x_i\partial_{x_i}$, the second
fundamental form is given by $I\!I_{\widetilde\nu}={\rm Id}$ and $\widetilde{H}=\tilde\nu$.
Therefore, Inequality \eqref{eq:inequalityimmersion} becomes
\begin{eqnarray*}
p q_{1,p-1} \int_{\partial M}(|\nu\lrcorner
d\omega|^2+|\iota^*(\delta\omega)|^2)d\mu_g&\leq&\\
&&\hspace{-3cm}\sum_{i=1}^{n+1}\int_M\Big|(2p-2n)\partial_{x_i}
^T\lrcorner\omega+2\langle\partial_ {
x_i}^\perp,\widetilde{\nu}\rangle\delta\omega-n
d(\langle\widetilde\nu,\partial_{x_i}
^\perp\rangle)\lrcorner\omega+\partial_
{ x_i } ^T\lrcorner\Delta\omega\Big|^2d\mu_g.
 \end{eqnarray*}
Now, an elementary computation shows that $d(\langle\widetilde\nu,\partial_{x_i}^\perp\rangle)=-\partial_{x_i}^T$. 
Therefore, the above inequality reduces to
\begin{eqnarray}\label{eq:inequalitysphere}
p q_{1,p-1} \int_{\partial M}(|\nu\lrcorner d\omega|^2+|\iota^*(\delta\omega)|^2)d\mu_g&\leq&\sum_{i=1}^{n+1}\int_M|(2p-n)\partial_{x_i}^T\lrcorner\omega+2\langle\partial_{x_i}^\perp,\widetilde{\nu}\rangle\delta\omega+\partial_{x_i}^T\lrcorner\Delta\omega|^2d\mu_g\nonumber\\
&=&\int_M\Big((2p-n)^2p|\omega|^2+4|\delta\omega|^2+p|\Delta\omega|^2\nonumber\\
&&\phantom{\int_M\Big(}
+2p(2p-n)
\langle\omega,\Delta\omega\rangle\Big)d\mu_g\nonumber\\
&=&\int_M(4|\delta\omega|^2+p|\Delta\omega+(2p-n)\omega|^2)d\mu_g.
\end{eqnarray}
Here, we use Identity \eqref{eq:scalinnerpro} and that
$$\sum_{i=1}^{n+1}\langle\partial_{x_i}^\perp,\widetilde{\nu}\rangle\langle\partial_{x_i}^T\lrcorner\omega,\delta\omega\rangle=\sum_{i=1}^{n+1}\langle\partial_{x_i},\widetilde{\nu}\rangle\langle\partial_{x_i}^T\lrcorner\omega,\delta\omega\rangle=\langle\widetilde{\nu}^T\lrcorner\omega,\delta\omega\rangle=0.$$
The same argument applies to
$\sum_{i=1}^{n+1}\langle\partial_{x_i}^\perp,\widetilde{\nu}
\rangle\langle\delta\omega,\partial_{x_i}^T\lrcorner\Delta\omega\rangle=0$.
Since Inequality \eqref{eq:inequalitysphere} is true for any $p$-eigenform
$\omega$, we apply it to the $(n-p)$-eigenform $*_M\omega$ to get
\begin{equation}\label{eq:33modified}
(n-p)q_{1,n-p-1} \int_{\partial M}(|\iota^*(\delta\omega)|^2+|\nu\lrcorner d\omega|^2)d\mu_g\leq \int_M(4|d\omega|^2+(n-p)|\Delta\omega+(2n-2p-n)\omega|^2)d\mu_g.
\end{equation}
Summing inequalities \eqref{eq:inequalitysphere} and \eqref{eq:33modified} and using the fact that $q_{1,n-p-1}=q_{1,p+1}$ yield the following:
\begin{eqnarray}\label{eq:inesphere}
\nonumber&&(p q_{1,p-1}+ (n-p) q_{1,p+1})\int_{\partial M}(|\nu\lrcorner
d\omega|^2+|\iota^*(\delta\omega)|^2)d\mu_g\\
&\leq&\int_M
(4|d\omega|^2+4|\delta\omega|^2+p|\Delta\omega+(2p-n)\omega|^2\nonumber+(n-p
)|\Delta\omega+(n-2p)\omega|^2)d\mu_g\nonumber\\
&=&\int_M
(4\langle\Delta\omega,
\omega\rangle+p|\Delta\omega+(2p-n)\omega|^2+(n-p)|\Delta\omega+(n-2p)\omega|^2)
d\mu_g.
\end{eqnarray}
Here, we used the fact that $\int_M\langle\Delta\omega,\omega\rangle
d\mu_g=\int_M (|d\omega|^2+|\delta\omega|^2)d\mu_g$ as a consequence of the
boundary condition $\omega_{|_{\partial M}}=0$. Now, the Bochner formula
applied to the eigenform $\omega$, with the help of the pointwise
inequality $|\nabla\alpha|^2\geq
\frac{1}{p+1}|d\alpha|^2+\frac{1}{n-p+1}|\delta\alpha|^2$ which is true for any
$p$-form $\alpha$ \cite[Lemme 6.8]{GM}, gives that
\begin{eqnarray*}
\int_M\langle\Delta\omega,\omega\rangle d\mu_g&=&\int_M|\nabla\omega|^2d\mu_g+\frac{1}{2}\int_M\Delta(|\omega|^2)d\mu_g+\int_M \langle\mathfrak{B}^{[p]}\omega,\omega\rangle d\mu_g\\
&\geq &\int_M\frac{1}{a(p,n)}(|d\omega|^2+|\delta\omega|^2)d\mu_g+\int_{\partial{M}}\langle\nabla_\nu\omega,\omega\rangle d\mu_g+\int_M p(n-p)|\omega|^2 d\mu_g\\
&=&\frac{1}{a(p,n)}\int_M\langle\Delta\omega,\omega\rangle d\mu_g+\int_M p(n-p)|\omega|^2 d\mu_g.
\end{eqnarray*}
Here, we have set $a(p,n)={\rm max}(p+1,n-p+1)$ and used the fact that the Bochner operator $\mathfrak{B}^{[p]}$ is equal to $p(n-p)\mathrm{Id}$ on the round sphere $\mathbb{S}^n$, see e.g. \cite[Cor. 2.6]{GM} and \cite[Rem. 6.15]{GM}. Thus, we deduce that
$$p(n-p)\|\omega\|^2_{L^2(M)}\leq \left(1-\frac{1}{a(p,n)}\right)\int_M\langle
\Delta\omega,\omega\rangle d\mu_g\leq
\left(1-\frac{1}{a(p,n)}\right)\|\Delta\omega\|_{L^2(M)}\|\omega\|_{L^2(M)},$$
so that
\begin{equation}\label{eq:caushw}
\|\omega\|_{L^2(M)}\leq \frac{a(p,n)-1}{p(n-p)a(p,n)}\|\Delta\omega\|_{L^2(M)}.
\end{equation}
Coming back to Inequality \eqref{eq:inesphere} and using again the
Cauchy-Schwarz inequality as well as the estimate
$\|\alpha+\beta\|^2_{L^2(M)}\leq 2(\|\alpha\|^2_{L^2(M)}+\|\beta\|^2_{L^2(M)})$, we obtain
\begin{eqnarray*}
(p q_{1,p-1}+ (n-p) q_{1,p+1})\int_{\partial M}(|\nu\lrcorner
d\omega|^2+|\iota^*(\delta\omega)|^2)d\mu_g&\leq &
4\|\Delta\omega\|_{L^2(M)}\|\omega\|_{L^2(M)}+2p\|\Delta\omega\|^2_{L^2(M)}
\\&&\hspace{-1.5cm}+2p(2p-n)^2\|\omega\|_{L^2(M)}^2\nonumber+2(n-p)\|\Delta\omega\|^2_{L^2(M)}
\\&&\hspace{-1.5cm}+2(n-p)(n-2p)^2\|\omega\|^2_{L^2(M)}\\
&\bui{\leq}{\eqref{eq:caushw}}& C(p,n)\|\Delta\omega\|^2_{L^2(M)},
\end{eqnarray*}
where $C(p,n)$ is the constant given by
\begin{eqnarray}\label{eq:c(p,n)}
C(p,n)&=&4\frac{a(p,n)-1}{p(n-p)a(p,n)}+2n+(2p(2p-n)^2+2(n-p)(n-2p)^2))(\frac{a(p,n)-1}{p(n-p)a(p,n)})^2\nonumber\\
&=& 4\frac{a(p,n)-1}{p(n-p)a(p,n)}+2n+2n(2p-n)^2(\frac{a(p,n)-1}{p(n-p)a(p,n)})^2.
\end{eqnarray}
Finally Characterization \eqref{eq:char1} allows to get the result. 
Notice that, if \eqref{eq:ineqq1psphere} were an equality, then by the limiting case in the Cauchy-Schwarz inequality, the form $\Delta\omega$ would be parallel to $\omega$. 
But by $\Delta^2\omega=0$, this would imply that $\Delta\omega=0$. 
Because of $\omega_{|_{\partial M}}=0$, we would deduce from \cite{Anne89} that $\omega=0$ on $M$, which is a contradiction.
Therefore, \eqref{eq:ineqq1psphere} always remains strict.
\hfill$\square$

\section{Robin vs. Dirichlet and Neumann eigenvalue problems}\label{sec:robinproblem}
In this section, we will establish estimates for the Robin eigenvalue problem on differential forms defined in \cite{EGH}. We mainly generalize some results in \cite{Ko} to differential forms. 
For this purpose, we recall the Robin problem on forms. Let $(M^n,g)$ be a compact Riemannian manifold with smooth boundary $\partial M$.
Fix a positive real number $\tau$.
Then the boundary value problem \cite{EGH}
\begin{equation}\label{Robin-forms}
\left\{\begin{array}{lll}\Delta \omega&=\lambda\omega &\textrm{ on } M\\ \iota^*(\nu\lrcorner\,d\omega-\tau\omega)&=0&\textrm{ on }\partial M\\
\nu\lrcorner\,\omega&=0&\textrm{ on }\partial M\end{array}\right.
\end{equation}
 is elliptic  and self-adjoint. It admits an increasing unbounded sequence of positive real eigenvalues with finite multiplicities $$\lambda_{1,p}(\tau)\le \lambda_{2,p}(\tau)\leq\cdots$$
The first eigenvalue $\lambda_{1,p}(\tau)$ of the Robin boundary problem {\rm(\ref{Robin-forms})} can be characterized as follows:
\begin{equation}\label{eq:cararobin}
\lambda_{1,p}(\tau)={\rm inf}\left\{\frac{\displaystyle\int_M \left(|d\omega|^2 +|\delta \omega|^2\right)d\mu_g+\tau\int_{\partial M} |\iota^*\omega|^2 d\mu_g}{\displaystyle\int_M|\omega|^2 d\mu_g}\right\},
\end{equation}
where $\omega$ runs over all non-identically vanishing $p$-forms on $M$ such that $\nu\lrcorner\omega=0$. When the parameter $\tau$ tends to $0$, the Robin problem \eqref{Robin-forms} reduces to the Neumann boundary problem, that is
\begin{equation}\label{Neumann-forms}
\left\{\begin{array}{lll}\Delta \omega&=\lambda\omega &\textrm{ on } M\\ \nu\lrcorner\,d\omega&=0&\textrm{ on }\partial M\\
\nu\lrcorner\,\omega&=0&\textrm{ on }\partial M.\end{array}\right.
\end{equation}
Notice that the first eigenvalue $\lambda_{1,p}^N$ of \eqref{Neumann-forms} is nonnegative and the kernel of the operator \eqref{Neumann-forms} is isomorphic to the so-called absolute de Rham cohomology, which is defined by
$$H_A^p(M)=\{\omega\in \Omega^p(M)|\,\, d\omega=\delta\omega=0\,\, \text{on}\,\, M\,\,\, \text{and}\,\,\, \nu\lrcorner\omega=0\,\, \text{on}\,\, \partial M\}.$$
When $\tau\to \infty$, the Robin problem \eqref{Robin-forms} reduces to the Dirichlet boundary problem
\begin{equation}\label{Dirichlet-forms}
\left\{\begin{array}{lll}\Delta \omega&=\lambda\omega &\textrm{ on } M\\ \omega&=0&\textrm{ on }\partial M.\\
\end{array}\right.
\end{equation}
By \cite{Anne89}, the first eigenvalue $\lambda_{1,p}^D$ of problem \eqref{Dirichlet-forms} is positive.
We also have the estimate \cite[Prop. 5.4]{EGH}
$$\lambda_{1,p}^N\leq \lambda_{1,p}(\tau)\leq \lambda_{1,p}^D.$$
In the following, we will establish another estimates for $\lambda_{1,p}(\tau)$ in terms of the Neumann and Dirichlet ones. 
We have:

\begin{thm} \label{est_Robin_Dir_Neu} Let $(M^n,g)$ be a compact Riemannian manifold with smooth boundary. We have the following estimates for the first eigenvalue of the Robin boundary problem:
\begin{enumerate}
\item Assume that the
absolute de Rham cohomology $H_A^p(M)$ does not vanish. We denote by $\omega_D$ an eigenform of the Dirichlet boundary problem associated to $\lambda^D_{1,p}$ and let $\omega_0$ be the orthogonal projection of $\omega_D$ on the space $H_A^p(M)$, assumed to be nonzero. Then
$$\dfrac{1}{\lambda_{1,p}(\tau)} \ge
\dfrac{1}{\lambda^D_{1,p}}+\dfrac{\|\omega_0\|^4_{L^2(M)}}{\tau\|\omega_0\|^2_{
L^2(\partial M)}} .$$
\item Assume that the first eigenvalue $\lambda^N_{1,p}$ of the Neumann boundary problem is positive, then
$$\dfrac{1}{\lambda_{1,p}(\tau)} \ge \dfrac{1}{\lambda^N_{1,p}}-\dfrac{\tau\alpha_N (\lambda^D_{1,p}-\lambda^N_{1,p})}{\lambda^N_{1,p}(\tau \alpha_N \lambda^D_{1,p}+\lambda^N_{1,p}(\lambda^D_{1,p}-\lambda^N_{1,p}))},$$
where $\alpha_N = \dfrac{\|\omega_N\|^2_{L^2(\partial
M)}}{\|\omega_N\|^2_{L^2(M)}}$ and $\omega_N$ is being an eigenform of the
Neumann boundary problem associated to $\lambda^N_{1,p}$.
\end{enumerate}
\end{thm}

\noindent{\it Proof.}
We begin with the proof of the first point.
Let $\omega_D$ be a $p$-eigenform associated to the first eigenvalue $\lambda_{1,p}^D$ that is assumed to be of $L^2$-norm equal to $1$. 
Let $\omega_0$ be the orthogonal projection of $\omega_D$ to $H_A^p(M)$. For any real number $t$, we consider the $p$-form
  $$\omega_t =\omega_D + t \, \omega_0.$$
Clearly, we have that $\nu\lrcorner\,\omega_t=0$. Therefore, by the characterization of the first eigenvalue $\lambda_{1,p}(\tau)$, we have that
   $$\lambda_{1,p}(\tau)\le \dfrac{\displaystyle\int_M \left(|d\omega_t|^2 +|\delta \omega_t|^2\right)d\mu_g+\tau \int_{\partial M}  |\omega_t|^2d\mu_g}{\displaystyle\int_M|\omega_t|^2d\mu_g}.$$
By the definition of the form $\omega_0$, we have that
$$\|d\omega_t\|^2_{L^2(M)}+\|\delta\omega_t\|^2_{L^2(M)}=\lambda_{1,p}^D,\,
\|\omega_t\|^2_{L^2(\partial M)}=t^2\|\omega_0\|^2_{L^2(\partial M)}.$$
Also, we have that
$$ \|\omega_t\|^2_{L^2(M)}=1+t^2 \|\omega_0\|^2_{L^2(M)}+2 t
\|\omega_0\|^2_{L^2(M)}.$$
The last term comes from the fact that $\omega_0$ is the orthogonal projection of $\omega_D$. Thus by plugging in the above inequality, we get that
   $$\lambda_{1,p}(\tau)\le \dfrac{\lambda^D_{1,p} + t^2 \tau
\|\omega_0\|^2_{L^2(\partial M)}}{1+t^2 \|\omega_0\|^2_{L^2(M)}+2 t
\|\omega_0\|^2_{L^2(M)}}.$$
Now, we take the inverse of this last inequality, then add and subtract the
term $\dfrac{t^2\tau  \|\omega_0\|^2_{L^2(\partial M)}}{\lambda^D_{1,p}}$ in
the numerator to find that
$$\dfrac{1}{\lambda_{1,p}(\tau)} \ge \dfrac{1}{\lambda^D_{1,p}} + \dfrac{t^2
\left(\|\omega_0\|^2_{L^2(M)} -\dfrac{\tau  \|\omega_0\|^2_{L^2(\partial
M)}}{\lambda^D_{1,p}} \right)+2 t  \|\omega_0\|^2_{L^2(M)}}{\lambda^D_{1,p} +
t^2  \tau  \|\omega_0\|^2_{L^2(\partial M)}}. $$
Since this is true for any real number $t$, then we deduce that $\dfrac{1}{\lambda_{1,p}(\tau)} \ge \dfrac{1}{\lambda^D_{1,p}}+\mathop{\rm sup}\limits_{\R}(f)$ where $f$ is the function given by
 $$f(t)= \dfrac{t^2 \left(\|\omega_0\|^2_{L^2(M)} -\dfrac{\tau
\|\omega_0\|^2_{L^2(\partial M)}}{\lambda^D_{1,p}} \right)+2 t
\|\omega_0\|^2_{L^2(M)}  }{\lambda^D_{1,p} + t^2  \tau
\|\omega_0\|^2_{L^2(\partial M)}}=\dfrac{At^2+Bt}{Ct^2+D}$$
with
 \begin{eqnarray*}
   A = \|\omega_0\|^2_{L^2(M)} -\dfrac{\tau  \|\omega_0\|^2_{L^2(\partial
M)}}{\lambda^D_{1,p}}, && B=2\|\omega_0\|^2_{L^2(M)},  \\
   C=  \tau  \|\omega_0\|^2_{L^2(\partial M)}, && D= \lambda^D_{1,p}.
 \end{eqnarray*}
It is easy to check that the supremum of $f$ is attained at
$t_0=\dfrac{AD +\sqrt{A^2D^2+B^2CD}}{BC}$ which is equal to
$$\sup_\R f= f(t_0)=\dfrac{t_0(At_0+B)}{Ct_0^2+D}=\dfrac{t_0(At_0+B)}{\dfrac{2ADt_0}{B}+D+D}=\dfrac{Bt_0}{2D}=\dfrac{A +\sqrt{A^2+\dfrac{B^2C}{D}}}{2C}.$$
Now, by replacing $A,B,C$ and $D$ by their values, we estimate
\begin{eqnarray*}
  A^2+\dfrac{B^2C}{D} &=& \left(\|\omega_0\|^2_{L^2(M)} -\dfrac{\tau
\|\omega_0\|^2_{L^2(\partial M)}}{\lambda^D_{1,p}}
\right)^2+\dfrac{4\|\omega_0\|^4_{L^2(M)} \tau  \|\omega_0\|^2_{L^2(\partial
M)}}{\lambda^D_{1,p}}\\
   &=& \left(\|\omega_0\|^2_{L^2(M)} -\dfrac{\tau  \|\omega_0\|^2_{L^2(\partial M)}}{\lambda^D_{1,p}} -2 \|\omega_0\|^4_{L^2(M)}  \right)^2-4\|\omega_0\|^8_{L^2(M)}+4\|\omega_0\|^6_{L^2(M)} \\
   &\geq& \left(\|\omega_0\|^2_{L^2(M)} -\dfrac{\tau  \|\omega_0\|^2_{L^2(\partial M)}}{\lambda^D_{1,p}} -2 \|\omega_0\|^4_{L^2(M)}  \right)^2,
\end{eqnarray*}
since $\|\omega_0\|^2_{L^2(M)} \leq \|\omega_D\|^2_{L^2(M)}=1$. Then
\begin{eqnarray*}
\sup_\R f &\geq&\dfrac{ \left(\|\omega_0\|^2_{L^2(M)} -\dfrac{\tau  \|\omega_0\|^2_{L^2(\partial M)}}{\lambda^D_{1,p}} \right) +\left| \,\|\omega_0\|^2_{L^2(M)} -\dfrac{\tau  \|\omega_0\|^2_{L^2(\partial M)}}{\lambda^D_{1,p}} -2 \|\omega_0\|^4_{L^2(M)} \,\right|}{2\tau  \|\omega_0\|^2_{L^2(\partial M)}}\\
&\geq& \dfrac{\|\omega_0\|^4_{L^2(M)}} {\tau  \|\omega_0\|^2_{L^2(\partial M)}}.
\end{eqnarray*}
This shows the required estimate. 
To prove the second inequality, let $\omega_D$ (resp. $\omega_N$) an eigenform of the Dirichlet (resp. Neumann) boundary problem associated to $\lambda^D_{1,p}$ (resp. $\lambda^N_{1,p}$) such that $\|\omega_D\|_{L^2(M)}=\|\omega_N\|_{L^2(M)}=1$. For any nonnegative number $s$, we consider the $p$-form $\omega_{s} =s\omega_D + \omega_N$, which clearly satisfies $\nu\lrcorner\,\omega_{s}=0$. 
In order to use the Rayleigh inequality for the eigenvalue $\lambda_{1,p}(\tau)$, we compute
\begin{eqnarray*}
\int_M (|d\omega_{s}|^2+|\delta\omega_s|^2) d\mu_g &=& s^2 \int_M (|d\omega_D|^2 + |\delta\omega_D|^2)d\mu_g+\int_M (|d\omega_N|^2+|\delta\omega_N|^2)d\mu_g\\&& +2s \int_M \langle d\omega_D, d \omega_N \rangle d\mu_g+ 2s\int_M \langle \delta\omega_D, \delta \omega_N \rangle d\mu_g \\
&=&s^2 \lambda^D_{1,p}+ \lambda^N_{1,p}+2 s \lambda^N_{1,p} (\omega_D, \omega_N)_{L^2(M)}.
\end{eqnarray*}
In the last equality, we use the Stokes formula. Also, we have
$$\|\omega_{s}\|_{L^2(\partial M)}^2= \|\omega_N\|_{L^2(\partial M)}^2\,\,\,\text{and}\,\,\, \|\omega_{s}\|_{L^2(M)}^2= s^2+1 +2 s  (\omega_D,\omega_N)_{L^2(M)}.$$
Therefore, after replacing we get that
   $$\lambda_{1,p}(\tau) \le \dfrac{s^2\lambda^D_{1,p}  +  \lambda^N_{1,p} +2 s \lambda^N_{1,p} (\omega_D, \omega_N)_{L^2(M)}+ \tau \|\omega_N\|_{L^2(\partial M)}^2}{s^2+ 1+2 s  (\omega_D,\omega_N)_{L^2(M)}}.$$
As we did before, we take the inverse of this last inequality, then add and subtract the term $\frac{s^2\lambda_{1,p}^D+\tau\|\omega\|^2_{L^2(\partial M)}}{\lambda^N_{1,p}}$ in the numerator to get that
$$\dfrac{1}{\lambda_{1,p}(\tau)} \ge  \dfrac{1}{\lambda^N_{1,p}}
   + \dfrac{s^2 \left( 1-\dfrac{\lambda^D_{1,p}}{\lambda^N_{1,p}}  \right) -\dfrac{\tau \alpha_N}{\lambda^N_{1,p}}   }{s^2\lambda^D_{1,p}  + \lambda^N_{1,p} +2s \lambda^N_{1,p} + \tau \alpha_N}.$$
Here, we also use the fact that $|(\omega_D,\omega_N)_{L^2(M)}|\leq 1$ by the Cauchy-Schwarz inequality and the fact that $s\geq 0$.  In order to get the lower bound, we need to compute the supremum of the function  $g$ which is given by
$$g(s)=\dfrac{As^2+B}{Cs^2+Ds+E},$$
    with
    \begin{eqnarray*}
      A= 1-\dfrac{\lambda^D_{1,p}}{\lambda^N_{1,p}}, & B= -\dfrac{\tau \alpha_N}{\lambda^N_{1,p}} , &  \\
      C=\lambda^D_{1,p}, & D= 2 \lambda^N_{1,p}, & E= \lambda^N_{1,p} + \tau \alpha_N.
    \end{eqnarray*}
A direct computation shows that the supremum of the function $g$ is attained at $s_2= \dfrac{ \tau \alpha_N}{\lambda^D_{1,p}-\lambda^N_{1,p}}$ and thus
$${\rm sup}_{\R}(g)=g(s_2)=-\dfrac{\tau \alpha_N(\lambda^D_{1,p}-\lambda^N_{1,p})}{\lambda^N_{1,p}(\tau \alpha_N \lambda^D_{1,p} +\lambda^N_{1,p}(\lambda^D_{1,p}-\lambda^N_{1,p}))}.$$
This finishes the proof of the theorem.
\hfill$\square$

In the following result, we will give a gap inequality between the eigenvalues of the Robin Laplacian under some curvature conditions. 
Let $(M^n,g)$ be a Riemannian manifold with smooth boundary and let $\eta_1(x), \eta_2(x),\cdots,\eta_{n-1}(x)$ be the principal curvatures at a point $x$ of the boundary (i.e. eigenvalues of the second fundamental form of the Weingarten tensor $S$). We assume that $\eta_1(x)\leq \eta_2(x)\leq\cdots\leq \eta_{n-1}(x)$. For any integer $p\in \{1,\cdots,n-1\}$, the $p$-curvatures $\sigma_p(x)$ are defined as $\sigma_p(x)=\eta_1(x)+\cdots+\eta_p(x)$. One can easily check that for any integer $p$ and $q$ with $p\leq q$, we have that $\frac{\sigma_p(x)}{p}\leq \frac{\sigma_q(x)}{q}$ with equality if and only if $\eta_1(x)=\eta_2(x)=\cdots=\eta_q(x)$. Hence, we deduce that $H\geq \frac{\sigma_p(x)}{p}$ for any integer $p\in \{1,\cdots,n-1\}$. In the next theorem, we set
$$\sigma_p=\mathop{\rm inf}\limits_{x\in \partial M}\sigma_p(x)$$
We state the result which generalizes \cite[Thm. 5.8]{EGH}.

\begin{thm}\label{t:gabRobin} Let $M$ be a compact domain in $\R^n$. Fix an integer number $q\in \{1,\cdots,n-1\}$ and Let $\omega$ be a $q$-eigenform of the Robin Laplacian. If $\sigma_p>0$ for some $p\leq q$, then we have
$$\frac{\|\omega\|^2_{L^2(\partial M)}}{\|\omega\|^2_{L^2(M)}}\leq \frac{1}{\sigma_p}(\lambda_{1,q}-\lambda_{1,q-p}).$$
\end{thm}

{\it Proof.} We mainly follow the same computations as in \cite[Thm. 5.8]{EGH}. Let $\omega$ be a $q$-eigenform of the Robin Laplacian and, for any $p\leq q$, consider the $(q-p)$-form $\phi_{i_1,\cdots,i_p}:=\partial_{x_{i_1}}\lrcorner\cdots\lrcorner \partial_{x_{i_p}}\lrcorner\omega$, for $i_k=1,\cdots,n$ with $k=1,\cdots,p$. Clearly, we have that $\nu\lrcorner\phi_{i_1,\cdots,i_p}=0$. 
Hence by the characterization \eqref{eq:cararobin} of the first eigenvalue, we get that
\begin{equation}\label{eq:lamdapq}
\lambda_{1,q-p}(\tau)\int_M|\phi_{i_1,\cdots,i_p}|^2d\mu_g\leq \int_M(|d\phi_{i_1,\cdots,i_p}|^2+|\delta\phi_{i_1,\cdots,i_p}|^2)d\mu_g+\tau\int_{\partial M}|\phi_{i_1,\cdots,i_p}|^2d\mu_g.
\end{equation}
Next we sum over $i_1,\cdots, i_p$. 
We begin with the l.h.s. 
Applying successively \eqref{eq:scalinnerpro}, we have
\begin{eqnarray} \label{eq:normphi}
\sum_{i_1,\cdots,i_p}|\phi_{i_1,\cdots,i_p}|^2&=&\sum_{i_1,\cdots,i_p}|\partial_{x_{i_1}}\lrcorner\cdots\lrcorner \partial_{x_{i_p}}\lrcorner\omega|^2\nonumber\\
&=&(q-(p-1))\cdots q\cdot|\omega|^2=\frac{q!}{(q-p)!}|\omega|^2.
\end{eqnarray}
For the r.h.s., we first compute
\begin{eqnarray}\label{eq:deltphi}
\sum_{i_1,\cdots,i_p}|\delta\phi_{i_1,\cdots,i_p}|^2&=&\sum_{i_1,\cdots,i_p}|\partial_{x_{i_1}}\lrcorner\cdots\lrcorner \partial_{x_{i_p}}\lrcorner\delta\omega|^2\bui{=}{\eqref{eq:scalinnerpro}}\frac{(q-1)!}{(q-1-p)!}|\delta\omega|^2.
\end{eqnarray}
In order to compute the term $\sum_{i_1,\cdots,i_p}|d\phi_{i_1,\cdots,i_p}|^2$, we proceed as in Equation \eqref{eq:splitLie}. 
Using repeatedly the identity $d(X\lrcorner\alpha)=\nabla_X\alpha-X\lrcorner d\alpha$, true for any parallel vector field $X$, we get
\begin{eqnarray*}
d\phi_{i_1,\cdots,i_p}=(-1)^p\partial_{x_{i_1}}\lrcorner\cdots\lrcorner \partial_{x_{i_p}}\lrcorner d\omega+\sum_{l=1}^p(-1)^{l+1}\partial_{x_{i_1}}\lrcorner\cdots\lrcorner\widehat{\partial_{x_{i_l}}}\lrcorner\cdots \lrcorner\partial_{x_{i_p}}\lrcorner\nabla_{\partial_{x_{i_l}}}\omega.
\end{eqnarray*}
Thus, we find that the sum $\sum_{i_1,\cdots,i_p}|d\phi_{i_1,\cdots,i_p}|^2$ is equal to
\begin{eqnarray*}
&&\sum_{i_1,\cdots,i_p}|\partial_{x_{i_1}}\lrcorner\cdots\lrcorner \partial_{x_{i_p}}\lrcorner d\omega|^2+\sum_{i_1,\cdots,i_p}|\sum_{l=1}^p(-1)^{l+1}\partial_{x_{i_1}}\lrcorner\cdots\lrcorner\widehat{\partial_{x_{i_l}}}\lrcorner\cdots \lrcorner\partial_{x_{i_p}}\lrcorner\nabla_{\partial_{x_{i_l}}}\omega|^2\\
&&+2(-1)^p\sum_{i_1,\cdots,i_p}\sum_{l=1}^p(-1)^{l+1}\langle\partial_{x_{i_1}}\lrcorner\cdots\lrcorner \partial_{x_{i_p}}\lrcorner d\omega,\partial_{x_{i_1}}\lrcorner\cdots\lrcorner\widehat{\partial_{x_{i_l}}}\lrcorner\cdots \lrcorner\partial_{x_{i_p}}\lrcorner\nabla_{\partial_{x_{i_l}}}\omega\rangle\\
&\bui{=}{\eqref{eq:scalinnerpro}}&\frac{(q+1)!}{(q+1-p)!}|d\omega|^2+\sum_{i_1,\cdots,i_p}\sum_{l}\langle\partial_{x_{i_1}}\lrcorner\cdots\lrcorner\widehat{\partial_{x_{i_l}}}\lrcorner\cdots \lrcorner\partial_{x_{i_p}}\lrcorner\nabla_{\partial_{x_{i_l}}}\omega,\partial_{x_{i_1}}\lrcorner\cdots\lrcorner\widehat{\partial_{x_{i_l}}}\lrcorner\cdots \lrcorner\partial_{x_{i_p}}\lrcorner\nabla_{\partial_{x_{i_l}}}\omega\rangle\\&&-2\sum_{i_1,\cdots,i_p}\sum_{l<s}(-1)^{l+s}\langle\partial_{x_{i_1}}\lrcorner\cdots\lrcorner\widehat{\partial_{x_{i_l}}}\lrcorner\cdots \lrcorner\partial_{x_{i_p}}\lrcorner\nabla_{\partial_{x_{i_l}}}\omega,\partial_{x_{i_1}}\lrcorner\cdots\lrcorner\widehat{\partial_{x_{i_s}}}\lrcorner\cdots \lrcorner\partial_{x_{i_p}}\lrcorner\nabla_{\partial_{x_{i_s}}}\omega\rangle\\
&&-2\sum_{i_1,\cdots,i_p}\sum_{l=1}^p\langle\partial_{x_{i_1}}\lrcorner\cdots\widehat{\partial_{x_{i_l}}}\lrcorner\cdots \lrcorner\partial_{x_{i_p}}\lrcorner (\partial_{x_{i_l}}\lrcorner d\omega),\partial_{x_{i_1}}\lrcorner\cdots\lrcorner\widehat{\partial_{x_{i_l}}}\lrcorner\cdots \lrcorner\partial_{x_{i_p}}\lrcorner\nabla_{\partial_{x_{i_l}}}\omega\rangle\\
&=&\frac{(q+1)!}{(q+1-p)!}|d\omega|^2+\frac{pq!}{(q-p+1)!}|\nabla\omega|^2\\&&-2\sum_{i_1,\cdots,i_p}\sum_{l<s}(q-1)-(p-3))\cdots(q-1)\langle\partial_{x_{i_s}}\lrcorner\nabla_{\partial_{x_{i_l}}}\omega,\partial_{x_{i_l}}\lrcorner\nabla_{\partial_{x_{i_s}}}\omega\rangle-2\frac{pq!}{(q-p+1)!}|d\omega|^2.
\end{eqnarray*}
In this lenghty computation, we used the fact that $\sum_i\langle\partial_{x_i}\lrcorner\alpha,\partial_{x_i}\lrcorner\beta\rangle=p\langle\alpha,\beta\rangle$ for any $p$-forms $\alpha$ and $\beta$. We also make use of the formula $d=\sum_i\partial_{x_i}\wedge\nabla_{\partial_{x_i}}$. Now, one can easily check by using the expression of $d$ and $\nabla$ that the sum term in the above equality is equal to
$$\sum_{i_1,\cdots,i_p}\sum_{l<s}\langle\partial_{x_{i_s}}\lrcorner\nabla_{\partial_{x_{i_l}}}\omega,\partial_{x_{i_l}}\lrcorner\nabla_{\partial_{x_{i_s}}}\omega\rangle=\begin{pmatrix}p\\2\end{pmatrix}(|\nabla\omega|^2-|d\omega|^2).$$
Hence, after simplifying, we find that
\begin{equation}\label{eq:dphi}
\sum_{i_1,\cdots,i_p}|d\phi_{i_1,\cdots,i_p}|^2=\frac{p(q-1)!}{(q-p)!}|\nabla\omega|^2+\frac{(q-1)!}{(q-p-1)!}|d\omega|^2.
\end{equation}
Plugging Equalities \eqref{eq:normphi}, \eqref{eq:deltphi} and \eqref{eq:dphi} into Inequality \eqref{eq:lamdapq} yields after simplifying by $\frac{(q-1)!}{(q-p)!}$
\begin{eqnarray}\label{eq:lambda1pq}
\lambda_{1,q-p}(\tau)q\int_M|\omega|^2 d\mu_g&\leq& p\int_M|\nabla\omega|^2d\mu_g+(q-p)\int_M(|d\omega|^2+|\delta\omega|^2)d\mu_g\nonumber\\&&+\tau q\int_{\partial M} |\omega|^2d\mu_g.
\end{eqnarray}
Now the Bochner formula $\Delta=\nabla^*\nabla+\mathfrak{B}^{[q]}$ applied to the form $\omega$ gives after integrating that
$$\lambda_{1,q}(\tau)\int_M|\omega|^2d\mu_g=\int_M\langle\Delta\omega,\omega\rangle d\mu_g=\int_M|\nabla\omega|^2 d\mu_g+\int_{\partial M}\langle\nabla_\nu\omega,\omega\rangle d\mu_g.$$
By using $\langle\nabla_\nu\omega,\omega\rangle=\langle S^{[q]}(\iota^*\omega),\iota^*\omega\rangle+\tau|\omega|^2$ as proved in \cite{GS}, we get that
$$\int_M|\nabla\omega|^2d\mu_g=\lambda_{1,q}(\tau)\int_M|\omega|^2d\mu_g-\int_{\partial M}(\langle S^{[q]}(\iota^*\omega),\iota^*\omega\rangle+\tau|\omega|^2)d\mu_g.$$
Plugging this last equality into Inequality \eqref{eq:lambda1pq}, we get after using the pointwise inequality $\langle S^{[q]}(\iota^*\omega),\iota^*\omega\rangle\geq \sigma_q|\omega|^2$ that
\begin{eqnarray*}
\lambda_{1,q-p}(\tau)q\int_M|\omega|^2 d\mu_g&\leq& p\lambda_{1,q}(\tau)\int_M |\omega|^2d\mu_g-p\int_{\partial M}(\sigma_q+\tau)|\omega|^2d\mu_g\nonumber\\&&+(q-p)\left(\lambda_{1,q}(\tau)\int_M|\omega|^2d\mu_g-\tau\int_{\partial M}|\omega|^2d\mu_g\right)+\tau q\int_{\partial M} |\omega|^2d\mu_g\\
&=&q\lambda_{1,q}(\tau)\int_M |\omega|^2d\mu_g-p\sigma_q\int_{\partial M}|\omega|^2d\mu_g\\
&\leq &q\lambda_{1,q}(\tau)\int_M |\omega|^2d\mu_g-q\sigma_p\int_{\partial M}|\omega|^2d\mu_g.
\end{eqnarray*}
In the last inequality, we used the fact that $\frac{\sigma_p}{p}\leq \frac{\sigma_q}{q}$. This finishes the proof of the theorem.
\hfill$\square$

\section{Robin and Steklov operators vs. biharmonic Steklov}\label{sec:RSvsS}

In this section, we relate the eigenvalues of the Robin problem to those of the biharmonic Steklov operator. 
This extends the result stated in \cite[Thm. 1.17]{Ko} for functions.

\begin{thm}\label{est_Robin_Dir_BS} Let $(M^n,g)$ be a compact Riemannian manifold with smooth boundary. We have for any $\tau>0$ the estimate
$$\dfrac{1}{\lambda_{1,p}(\tau)} \le \dfrac{1}{\lambda^D_{1,p}}+\dfrac{1}{\tau q_{1,p}},$$ where $\lambda^D_{1,p}$ is the first eigenvalue of the Dirichlet boundary problem \eqref{Dirichlet-forms}.
\end{thm}

\noindent{\it Proof.}  We mainly follow the computations done in \cite[Thm. 1.17]{Ko}.
Let $\omega$ be an eigenform of the Robin boundary problem {\rm(\ref{Robin-forms})} associated to $\lambda_{1,p}(\tau)$.
We denote by $\omega_1$ a solution to the problem
  \begin{equation}
\left\{\begin{array}{lll}\Delta \omega_1 &=0 &\textrm{ on } M\\ \iota^* \omega_1&=\iota^* \omega &\textrm{ on }\partial M\\
\nu\lrcorner\,\omega_1&=0&\textrm{ on }\partial M.\end{array}\right.
\end{equation}
Notice that such a problem admits a unique solution by \cite{Sc}. Now, let us consider  the $p$-form $\omega_2:=\omega-\omega_1$. It clearly satisfies
\begin{equation}
\left\{\begin{array}{lll}\Delta \omega_2 &=\Delta \omega &\textrm{ on } M\\ \omega_2&=0 &\textrm{ on }\partial M.
\end{array}\right.
\end{equation}
By using the triangle inequality, the characterization \eqref{eq:char2} and the one of the first eigenvalue $\lambda_{1,p}^D$ of the Dirichlet problem (corresponding to $\tau\to\infty$ in \eqref{eq:cararobin}), we have
\begin{eqnarray}
\label{normofomega}
  \|\omega\|_{L^2(M)} &\le& \|\omega_1\|_{L^2(M)}+\|\omega_2\|_{L^2(M)} \nonumber \\
   &\le & \sqrt{q_{1,p}^{-1}} \|\omega_1\|_{L^2(\partial M)} + \sqrt{(\lambda_{1,p}^D)^{-1}}\left(\|d \omega_2\|^2_{L^2(M)} + \|\delta \omega_2\|^2_{L^2(M)}     \right)^\frac{1}{2} \nonumber \\
   &\le & \sqrt{q_{1,p}^{-1}} \|\omega\|_{L^2(\partial M)} + \sqrt{(\lambda_{1,p}^D)^{-1}}\left(\|d \omega\|^2_{L^2(M)} + \|\delta \omega\|^2_{L^2(M)}     \right)^\frac{1}{2}.
\end{eqnarray}
Indeed, $\|\omega_1\|_{L^2(\partial M)}=\|\omega\|_{L^2(\partial M)}$ since $\omega_2=0$ on $\partial M$ and
\begin{eqnarray*}
  \|d \omega\|^2_{L^2(M)} + \|\delta \omega\|^2_{L^2(M)} &=& \|d \omega_1\|^2_{L^2(M)} +\|d \omega_2\|^2_{L^2(M)}+2(d \omega_1, d\omega_2)_{L^2(M)} \\
  && + \|\delta \omega_1\|^2_{L^2(M)}+ \|\delta \omega_2\|^2_{L^2(M)}+2(\delta \omega_1, \delta\omega_2)_{L^2(M)}\\
   &=& \|d \omega_1\|^2_{L^2(M)} + \|\delta \omega_1\|^2_{L^2(M)}+ \|d \omega_2\|^2_{L^2(M)}+ \|\delta \omega_2\|^2_{L^2(M)} \\
   && +2 \int_M \langle \delta d \omega_1, \omega_2 \rangle d\mu_g -2 \int_{\partial M} \langle  \nu\lrcorner\,\omega_1, \iota^* \omega_2 \rangle d\mu_g
  \\&& +2 \int_M \langle d\delta \omega_1, \omega_2 \rangle d\mu_g+2 \int_{\partial M} \langle  \iota^* \delta \omega_1,  \nu\lrcorner\, \omega_2 \rangle d\mu_g\\
   &=&  \|d \omega_1\|^2_{L^2(M)} + \|\delta \omega_1\|^2_{L^2(M)}+ \|d \omega_2\|^2_{L^2(M)}+ \|\delta \omega_2\|^2_{L^2(M)}  \\
   &\ge & \|d \omega_2\|^2_{L^2(M)}+ \|\delta \omega_2\|^2_{L^2(M)}.
\end{eqnarray*}
Now, we square both sides of Inequality \eqref{normofomega} to write
\begin{eqnarray*}
  \|\omega\|^2_{L^2(M)} &\le& q_{1,p}^{-1} \|\omega\|^2_{L^2(\partial M)}+ (\lambda_{1,p}^D)^{-1} \left(\|d \omega\|^2_{L^2(M)} + \|\delta \omega\|^2_{L^2(M)}     \right) \nonumber \\
  && + 2 \left( q_{1,p}^{-1} (\lambda_{1,p}^D)^{-1}  \|\omega\|^2_{L^2(\partial M)} \left(\|d \omega\|^2_{L^2(M)} + \|\delta \omega\|^2_{L^2(M)}     \right)       \right)^\frac{1}{2} \nonumber \\
  &\le& q_{1,p}^{-1} \|\omega\|^2_{L^2(\partial M)}+ (\lambda_{1,p}^D)^{-1} \left(\|d \omega\|^2_{L^2(M)} + \|\delta \omega\|^2_{L^2(M)}     \right) \nonumber \\
  && \tau^{-1} q_{1,p}^{-1}  \left(\|d \omega\|^2_{L^2(M)} + \|\delta \omega\|^2_{L^2(M)}     \right)
  + \tau (\lambda_{1,p}^D)^{-1}\|\omega\|^2_{L^2(\partial M)}
   \nonumber \\
  &= & (\tau^{-1} q_{1,p}^{-1} +  (\lambda_{1,p}^D)^{-1}) \left(\|d \omega\|^2_{L^2(M)} + \|\delta \omega\|^2_{L^2(M)} + \tau \|\omega\|^2_{L^2(\partial M)} \right).
  \end{eqnarray*}
In the second above inequality, we use the fact that $2\sqrt{ab}\leq \frac{a}{\tau}+\tau b$ for any real positive $\tau$. 
The characterization \eqref{eq:cararobin} allows to deduce the estimate.
\hfill$\square$

Now, we come back to the Serrin problem on forms.
We will use the existence of solution to this problem in harmonic domains carrying parallel forms to get an estimate for the eigenvalues of the absolute Dirichlet-to-Neumann operator introduced in \cite{RS} (see also \cite{K}). 
Recall the definition of this operator. Let $(M^n,g)$ be a compact Riemannian manifold with smooth boundary $\partial M$. 
Let $p\in\{0,\cdots, n-1\}$. 
Given any $p$-form $\omega$ on $\partial M$, there exists a unique $p$-form $\hat\omega$ on $M$ such that \cite{Sc}

\begin{equation}\label{dtn}
\left\{\begin{array}{lll}
\Delta \hat\omega&=0 &\textrm{ on } M\\
\iota^*\hat\omega &=\omega&\textrm{ on } \partial M\\
\nu\lrcorner\omega&=0&\textrm{ on }\partial M.\\
\end{array}\right.
\end{equation}

The form $\hat\omega$ is usually called the harmonic tangential extension of $\omega$. The Dirichlet-to-Neumann operator is then defined as $T^{[p]}:\Lambda^p(\partial M)\to \Lambda^p(\partial M)$, $\omega\mapsto -\nu\lrcorner d\hat\omega$.  
When $p=0$, this operator reduces to the classical Dirichlet-to-Neumann operator on functions, known as Steklov operator. It is shown in \cite{RS} that $T^{[p]}$ is an elliptic self-adjoint pseudo-differential operator with discrete spectrum consisting of eigenvalues
$$0\leq \nu_{1,p}(M)\leq \nu_{2,p}(M)\leq\cdots.$$
The kernel of this operator is isomorphic to the absolute de Rham cohomology $H_A^p(M)$ introduced in Section \ref{sec:robinproblem}. 
The dual problem to \eqref{dtn} (w.r.t. the Hodge star operator) is called the relative Dirichlet-to-Neumann operator and is defined by $T^{[p]}_D=(-1)^{p(n-1-p)}*_{\partial M}T^{[n-1-p]}*_{\partial M}$. If $\nu_{1,p}^D(M)$ is the first eigenvalue of $T^{[p]}_D$, we have
$$\nu_{1,p}^D(M)=\nu_{1,n-1-p}(M).$$
Also, we have the following characterization \cite{RaulotSavo2015} for the first eigenvalue $\nu_{1,p}^D(M)$: 
\begin{equation}\label{eq:carasteklov}
\nu_{1,p}^D(M)={\rm inf}\left\{\frac{\displaystyle\int_M(|d\phi|^2+|\delta\phi|^2)d\mu_g}{\displaystyle\int_{\partial M}|\phi|^2d\mu_g}|\,\, \phi\in\Omega^{p+1}(M), \iota^*\phi=0\right\}.
\end{equation}

\begin{thm} Let $(M^n,g)$ be a compact Riemannian manifold with smooth boundary. Assume that $M$ is a harmonic domain and carries a parallel $p$-form for some $p=1,\cdots, n-1$. 
If moreover $\sigma_p>0$ or $\sigma_{n-p}>0$, then
$${\rm min}(\nu_{1,p-1}(M),\nu_{1,n-1-p}(M))\leq \frac{{\rm Vol}(\partial M)}{{\rm Vol}(M)}.$$
\end{thm}

{\it Proof.} 
Let $\omega$ be any $p$-eigenform of the biharmonic Steklov operator associated to some eigenvalue, say $q$. We let the $(p+1)$-form $\phi:=d\omega$. Clearly, we have that $\iota^*\phi=0$. Hence, by Characterization \eqref{eq:carasteklov}, we get that
\begin{equation}\label{eq:inequalitynu1}
\nu_{1,p}^D(M)\displaystyle\int_{\partial M}|\nu\lrcorner d\omega|^2d\mu_g \leq \displaystyle\int_M|\delta d\omega|^2d\mu_g.
\end{equation}
Now, applying the same characterization \eqref{eq:carasteklov} to $d(*_M\omega)$, since $*_M\omega$ is also an eigenform of the biharmonic Steklov operator, yields the inequality
\begin{equation}\label{eq:inequalitynu2}
\nu_{1,n-p}^D(M)\displaystyle\int_{\partial M}|\iota^*\delta\omega|^2d\mu_g\leq \displaystyle\int_M|d\delta \omega|^2d\mu_g.
\end{equation}
Summing Inequalities \eqref{eq:inequalitynu1} and \eqref{eq:inequalitynu2} yields
$${\rm min}(\nu_{1,p}^D(M),\nu_{1,n-p}^D(M))\int_{\partial M}(|\nu\lrcorner d\omega|^2+|\iota^*\delta\omega|^2)d\mu_g\leq \int_M(|\delta d\omega|^2+|d\delta\omega|^2)d\mu_g.$$
Now a direct computation using the Stokes formula and the boundary conditions on $\omega$ gives that
$$\int_M|\Delta\omega|^2d\mu_g=\int_M(|d\delta\omega|^2+|\delta d\omega|^2)d\mu_g-\frac{2}{q^2}\int_{\partial M} \langle \nu\lrcorner \Delta\omega,\delta^{\partial M}(\iota^*\Delta\omega)\rangle d\mu_g.$$
Hence after plugging this equality into the above inequality, we deduce that
\begin{equation}\label{eq:inequalformespropres}
{\rm min}(\nu_{1,p}^D(M),\nu_{1,n-p}^D(M))\leq q+\frac{2}{q^2}\frac{\displaystyle\int_{\partial M} \langle \nu\lrcorner \Delta\omega,\delta^{\partial M}(\iota^*\Delta\omega)\rangle d\mu_g}{\displaystyle\int_{\partial M}(|\nu\lrcorner d\omega|^2+|\iota^*\delta\omega|^2)d\mu_g}.
\end{equation}
Notice that this inequality is true for any eigenform $\omega$ of the biharmonic Steklov operator. From Theorem \eqref{thm:harmonicdomain}, we know that when $M$ is a harmonic domain carrying a parallel $p$-form $\omega_0$, the form $\omega=f\cdot\omega_0$ is an eigenform associated to the eigenvalue $q=\frac{{\rm Vol}(\partial M)}{{\rm Vol}(M)}$. Hence, we will apply Inequality \eqref{eq:inequalformespropres} to the particular form $\omega=f\cdot\omega_0$. 
For this purpose, we will check the sign of the integral. 
Assume first that $\sigma_{n-p}>0$.
We estimate
\begin{eqnarray*}
\int_{\partial M} \langle \nu\lrcorner \Delta\omega,\delta^{\partial M}(\iota^*\Delta\omega)\rangle d\mu_g
&=&\int_{\partial M}\langle\nu\lrcorner\omega_0,\delta^{\partial M}(\iota^*\omega_0)\rangle d\mu_g\\
&=&\int_{\partial M}\langle\nu\lrcorner\omega_0,(S^{[p-1]}-(n-1)H)\nu\lrcorner\omega_0)\rangle d\mu_g\\
&\leq& \int_{\partial M} \left((\sigma_{n-1}-\sigma_{n-p})|\nu\lrcorner\omega_0|^2-\sigma_{n-1}|\nu\lrcorner\omega_0|^2\right) d\mu_g\\
&\leq &0.
\end{eqnarray*}
In the second equality, we use the identity \cite[Lem. 18]{RaulotSavo2011}
$$\delta^{\partial M}(\iota^*\omega_0)=\nu\lrcorner\nabla_\nu\omega_0+\iota^*(\delta\omega_0)+S^{[p-1]}(\nu\lrcorner\omega_0)-(n-1)H\nu\lrcorner\omega_0.$$
Also, we use the pointwise estimate $\langle S^{[p]}\alpha,\alpha\rangle\leq (\sigma_{n-1}-\sigma_{n-p-1})|\alpha|^2$ for any $p$-form $\alpha$. Hence, we deduce that
$${\rm min}(\nu_{1,p}^D(M),\nu_{1,n-p}^D(M))\leq \frac{{\rm Vol}(\partial M)}{{\rm Vol}(M)}.$$
Finally, the fact that $\nu_{1,p}^D(M)=\nu_{1,n-1-p}(M)$ and $\nu_{1,n-p}^D(M)=\nu_{1,p-1}(M)$ finishes the proof of the statement when $\sigma_{n-p}>0$.
If $\sigma_p>0$, then replacing $p$ by $n-p$ and by invariance of parallel forms under the Hodge star operator on $M$, we obtain the same inequality.
This concludes the proof.
\hfill$\square$

\begin{remark}
\rm
Notice that a similar estimate has been established in \cite[Cor. 15]{RS} for compact manifolds carrying parallel forms with the assumption $H_A^p(M)=H^p(M)=0$. The inequality is
$$\nu_{1,n-1-p}(M)+\nu_{1,p-1}(M)\leq \frac{{\rm Vol}(\partial M)}{{\rm Vol}(M)}.$$
\end{remark}
\section{Appendix}

\begin{lemma} Let $(M^n,g)$ be a compact Riemannian manifold with smooth boundary $\partial M$ and let $\nu$ be the inward unit normal vector field to the boundary. Consider the following boundary value problem
	\begin{equation}\label{eq:bvporder4}
	\left\{\begin{array}{lll}\Delta^2\omega&=f&\textrm{ on
	}M\\B_1 \omega&=\omega_1&\textrm{ on
	}\partial M\\B_2\omega&=\omega_2&\textrm{ on
	}\partial M\\B_3 \omega&=\omega_3&\textrm{ on
	}\partial M\end{array}\right.
	\end{equation}
	for given $f\in\Omega^p(M)$, $\omega_1\in\Gamma(\Lambda^pT^*M_{|_{\partial
			M}})$,
	$\omega_2\in\Omega^{p-1}(\partial M)$ and $\omega_3\in\Omega^p(\partial M)$ and where if $E_1:=\Lambda^pT^*M_{|_{\partial M}}$,
	$E_2:=\Lambda^{p-1}T^*\partial M$ and $E_3:=\Lambda^p T^*\partial M$, $ B_1\colon\Omega^p(M)\to\Gamma(E_1)$ such that $B_1\omega:=\omega_{|_{\partial M}}$,  $B_2\colon\Omega^p(M)\to\Gamma(E_2)$ and $B_3\colon\Omega^p(M)\to\Gamma(E_3)$ are either: \begin{enumerate}
		\item $B_2=\iota^*\delta\omega$ and $B_3=\nu\lrcorner d\omega$. In this case, $\eqref{eq:bvporder4}$ is elliptic in the
		sense of Lopatinski\u{\i}-Shapiro (see Definition $1.6.1$ in \cite{Sc}), self-adjoint and its kernel is reduced to $\{0\}$. Or
		\item $B_2\omega:=\nu\lrcorner\Delta\omega+q\iota^*(\delta\omega)$ and
		$B_3\omega:=\iota^*\Delta\omega-q\nu\lrcorner\,d\omega$ for some real constant
		$q$. 
In this case, problem $\eqref{eq:bvporder4}$ is elliptic in the sense of Lopatinski\u{\i}-Shapiro.
	\end{enumerate}
\end{lemma}
{\it Proof.}
Given any $v\in T_x^*\partial M\setminus\{0\}$ for a fixed $x\in\partial M$, we
consider the space
\[\mathcal{M}_{v}^+:=\left\{\textrm{bounded solutions }y=y(t)\textrm{ on
}\R_+\textrm{ to the ODE }\sigma_{\Delta^2}((-iv,\partial_t))y=0\right\}.\]
A direct computation shows that
\[\mathcal{M}_{v}^+=\mathrm{Span}\left(e^{-|v|t}(at+b)\cdot\omega_0\,|\,a,b\in\R
,\,\omega_0\in\Lambda^pT_x^*M_{|_{\partial M}}\right),\]
which is hence a space of dimension
$N:=2\left(\begin{array}{c}n\\p\end{array}\right)$.
We look at the pointwise map
\begin{eqnarray*}
	\mathcal{M}_{v}^+&\longrightarrow&\bigoplus_{j=1}^3 E_j\\
	y&\longmapsto&\left(\sigma_{B_1}((-iv,\partial_t))y,\sigma_{B_2}((-iv,
	\partial_t))y,\sigma_{B_3}((-iv,
	\partial_t))y\right)(0)
\end{eqnarray*}
which we want to show to be an isomorphism.
Note already that space dimensions are equal on both sides.
Since $\sigma_{B_1}((-iv,\partial_t))=\mathrm{Id}$,
$\sigma_{B_2}((-iv,
\partial_t))=-\partial_t\cdot\nu\lrcorner\cdot+iv\lrcorner\iota^*\cdot$ and
$\sigma_{B_3}((-iv,\partial
t))=\partial_t\cdot\iota^*+iv\wedge(\nu\lrcorner\cdot)$, we obtain that, for
any fixed $\omega_0\in\Lambda^pT_x^*M$, the element $e^{-|v|t}\cdot\omega_0$ of
$\mathcal{M}_{v}^+$ (corresponding to $a=0$ and $b=1$) is sent to
$(\omega_0,|v|\nu\lrcorner\omega_0+iv\lrcorner\iota^*\omega_0,
-|v|\iota^*\omega_0+iv\wedge(\nu\lrcorner\omega_0))$; and that the
element $te^{-|v|t}\cdot\omega_0$ of
$\mathcal{M}_{v}^+$ (corresponding to $a=1$ and $b=0$) is sent to
$(0,-\nu\lrcorner\omega_0,
\iota^*\omega_0)$.
Choosing a basis $(\omega_0^{(1)},\ldots,\omega_0^{(N)})$ of $\Lambda^pT_x^*M$,
the basis
$\left(e^{-|v|t}\cdot\omega_0^{(1)},\ldots,e^{-|v|t}\cdot\omega_0^{(N)},te^{
	-|v|t }\cdot\omega_0^{(1)},\ldots,te^{-|v|t}\cdot\omega_0^{(N)}\right)$ of
$\mathcal{M}_{v}^+$ will therefore be sent to a basis of $\bigoplus_{j=1}^3
E_j$.
This shows the map $\mathcal{M}_{v}^+\longrightarrow\bigoplus_{j=1}^3 E_j$ to
be an isomorphism.
Therefore \eqref{eq:bvporder4} is elliptic.\\
Using \eqref{eq:partialintDelta}, it is easy to see that \eqref{eq:bvporder4}
is also self-adjoint.
Moreover, the kernel of \eqref{eq:bvporder4} is reduced to $\{0\}$: namely, if
$\omega\in\Omega^p(M)$ solves \eqref{eq:bvporder4} with $f=0$ as well as
$\omega_1=\omega_2=\omega_3=0$, then \eqref{eq:partialintDeltacons} implies that
$\|\Delta\omega\|_{L^2(M)}=0$, from which $\omega=0$ on $M$ follows using
$\omega_{|_{\partial M}}=0$.\\
As a consequence, fixing $f=0$ as well as $\omega_1=0$, for any given
$(\omega_2,\omega_3)\in\Omega^{p-1}(\partial M)\oplus\Omega^p(\partial M)$,
there exists a unique $\omega\in\Omega^p(M)$
solving \eqref{eq:bvporder4}.
In particular, $\omega\in Z$.
This shows the map $Z\to
\Omega^p(\partial M)\oplus\Omega^{p-1}(\partial M)$,
$\omega\mapsto(\nu\lrcorner d\omega,\iota^*\delta\omega)$, to be onto. This proves $1.$\\
Changing the boundary operators $B_2$ and $B_3$ via
$B_2\omega:=\nu\lrcorner\Delta\omega+q\iota^*(\delta\omega)$ and
$B_3\omega:=\iota^*\Delta\omega-q\nu\lrcorner\,d\omega$ for some real contant
$q$ (which actually plays no role since it is only involved in the
first-order-terms of the b.c. and not in their principal symbols), we still get
elliptic boundary conditions for $\Delta^2$: for any $v\in T_x^*\partial
M\setminus\{0\}$,
the pointwise map
$\mathcal{M}_{v}^+\longrightarrow\bigoplus_{j=1}^3 E_j$ from above sends
$e^{-|v|t}\cdot\omega_0$ to $(\omega_0,0,0)$ and sends
$te^{-|v|t}\cdot\omega_0$ to
$(0,2|v|\nu\lrcorner\omega_0,2|v|\iota^*\omega_0)$, which shows that map to be
an isomorphism.
Therefore $B_1,B_2,B_3$ define elliptic boundary conditions for $\Delta^2$, this proves $2.$\hfill$\square$\\

\begin{lemma}\label{delta} Let $M^n\hookrightarrow \R^{n+m}$ be an isometric immersion and let $I\!I$ be the second fundamental form of the immersion. For all $X,Y\in TM$ and $N\in T^\perp M$, we have
\begin{equation}
\label{eq:delta}
(\nabla_X I\!I_N)(Y)=(\nabla_Y I\!I_N)(X)-I\!I_{\nabla_Y^\perp N}(X)+I\!I_{\nabla_X^\perp N}(Y),
\end{equation}
where $\nabla_X^\perp N:=(\nabla_X^{\R^{n+m}}N)^\perp$ defines the normal connection on $T^\perp M$.
As a consequence, by writing $\partial_{x_i}=\partial_{x_i}^T+\partial_{x_i}^\perp$ for all $i=1,\cdots,n+m$, the divergence of the endomorphism $I\!I_{\partial_{x_i}^{\perp}}$ is equal to
	\begin{equation}\label{Sec-term-delta}\delta(I\!I_{\partial_{x_i}^{\perp}})= -n d(\langle \widetilde{H},\partial_{x_i}^\perp\rangle)-n I\!I_{\widetilde{H}}(\partial_{x_i}^T)+\sum_{a=1}^m
	I\!I^2_{f_a}(\partial_{x_i}^T).	
	\end{equation}
Here $\{f_1,\cdots, f_m\}$ is a local orthonormal frame of $T^\perp M$ and $\widetilde{H}$ is the mean curvature field of the immersion.
\end{lemma}

{\it Proof.} Let $X,Y,Z$ be vector fields in $TM$  that we assume to be parallel at some point in $M$ and let $N \in T^\perp M$. We compute
\begin{eqnarray*}
(\nabla_X I\!I_N)(Y,Z)&=& X (I\!I_N(Y,Z)) \quad \\
&=&X \langle I\!I(Y,Z), N\rangle
\\
&=&  \langle \nabla_X^{\R^{n+m}}I\!I(Y,Z), N\rangle +  \langle I\!I(Y,Z),\nabla_X^{\R^{n+m}} N\rangle
\\
&=& \langle \nabla_Y^{\R^{n+m}} I\!I(X,Z), N\rangle +  \langle I\!I(Y,Z),\nabla_X^{\R^{n+m}} N\rangle\\
&=& Y(I\!I_N(X,Z))-\langle I\!I(X,Z),  \nabla_Y^{\R^{n+m}} N\rangle+\langle I\!I(Y,Z),\nabla_X^{\R^{n+m}} N\rangle\\
&=&(\nabla_Y I\!I_N)(X,Z)-I\!I_{\nabla_Y^\perp N}(X,Z)+I\!I_{\nabla_X^\perp N}(Y,Z).
\end{eqnarray*}
In the fourth equality, we use the Codazzi equation for submanifolds in $\R^{n+m}$. Hence we get Equality \eqref{eq:delta}. To find Equation \eqref{Sec-term-delta}, we  decompose $\partial_{x_i}=\partial_{x_i}^T+\partial_{x_i}^\perp$ for all $i=1,\cdots,n+m$. Then from the parallelism of the vector $\partial_{x_i}$ and the Gauss formula, we get that
$$\nabla^{\R^{n+m}}_X \partial_{x_i}^{\perp}=-\nabla^{\R^{n+m}}_X \partial_{x_i}^{T}= - \nabla_X \partial_{x_i}^{T} -I\!I(X, \partial_{x_i}^{T}),$$
where $\nabla$ is the Levi-Civita connection on $TM$. Thus, we deduce that
\begin{equation}\label{eq:nablaxiperp}
\nabla_X^\perp\partial_{x_i}^\perp=(\nabla^{\R^{n+m}}_X \partial_{x_i}^{\perp})^\perp=-I\!I(X, \partial_{x_i}^{T}).
\end{equation}
The divergence of the endomorphism $I\!I_{\partial_{x_i}^\perp}$ can be computed using Equation \eqref{eq:delta} with $N=\partial_{x_i}^\perp$. For $X\in TM$, we have
\begin{eqnarray*}
\delta(I\!I_{\partial_{x_i}^{\perp}})(X)&=& -\sum_{s=1}^n(\nabla_{e_s}I\!I_{\partial_{x_i}^{\perp}})(e_s,X)\\
&=&-\sum_{s=1}^n(\nabla_{e_s}I\!I_{\partial_{x_i}^{\perp}})(X,e_s)\\
&\bui{=}{\eqref{eq:delta}}&-\sum_s(\nabla_{X}I\!I_{\partial_{x_i}^{\perp}})(e_s,e_s)-\sum_{s=1}^n I\!I_{I\!I(X, \partial_{x_i}^{T})}(e_s,e_s)+ \sum_{s=1}^n  I\!I_{I\!I(e_s, \partial_{x_i}^{T})}(X,e_s) \nonumber\\
&=&-nX(\langle\widetilde H,\partial_{x_i}^\perp\rangle)-n I\!I_{\widetilde{H}}(X,\partial_{x_i}^T)+\sum_{a=1}^m\sum_{s=1}^n\langle I\!I(e_s,\partial_{x_i}^{T}),f_a\rangle \langle I\!I(X,e_s),f_a\rangle \\
&=&-nX(\langle\widetilde H,\partial_{x_i}^\perp\rangle)-n I\!I_{\widetilde{H}}(\partial_{x_i}^T,X)+ \sum_{a=1}^m
g(I\!I_{f_a}(\partial_{x_i}^T),I\!I_{f_a}(X)).
\end{eqnarray*}
This ends the proof of the lemma.
\hfill$\square$

\begin{prop} \label{prop:curvature} For any $p$-form $\omega$, the curvature term $\sum_{s=1}^n e_s \lrcorner R(e_s,X) \omega $ is equal to
$$\sum_{s=1}^n e_s \lrcorner R(e_s,X) \omega=- \sum_{a=1}^m I\!I_{f_a}(X)\lrcorner I\!I_{f_a}^{[p]}\omega+I\!I_{n\widetilde{H}}(X)\lrcorner\omega,$$ where $\{e_1,\cdots,e_n\}$ and $\{f_1,\cdots,f_m\}$ are respectively orthonormal basis of $TM$ and $T^\perp M$.
\end{prop}

{\it Proof.} In order to compute the curvature term, we use the Gau\ss{} equation. Indeed, for any $X,Y,Z,T\in TM$ we have
$$R(X,Y,Z,T)=-\langle I\!I(X,Z), I\!I(Y,T) \rangle + \langle I\!I(X,T), I\!I(Y,Z) \rangle,$$
which can be equivalently written  as
$$ R(X,Y)Z=-\sum_{a=1}^m g(I\!I_{f_a}(X),Z) I\!I_{f_a}(Y)+\sum_{a=1}^m g(I\!I_{f_a}(Y),Z) I\!I_{f_a}(X).$$
Now, due to linearity, we can consider that $\omega$ of the form $e_{i_1}\wedge\cdots\wedge e_{i_p}$. Then, we compute
\begin{eqnarray*}
&& \sum_{s=1}^n e_s \lrcorner R(e_s,X) \omega \\
&=& \sum_{s=1}^n \sum_{j=1}^p e_s \lrcorner (e_{i_1}\wedge \dots \wedge R(e_s,X)e_{i_j} \wedge \dots \wedge e_{i_p})\\
&=& -\sum_{a=1}^m\sum_{s=1}^n \sum_{j=1}^p  g(I\!I_{f_a}(e_s),e_{i_j}) e_s \lrcorner (e_{i_1}\wedge \dots \wedge \underbrace{I\!I_{f_a}(X)}_{j-th} \wedge \dots \wedge e_{i_p})\\
&&+\sum_{a=1}^m\sum_{s=1}^n \sum_{j=1}^p g(I\!I_{f_a}(X),e_{i_j}) e_s \lrcorner (e_{i_1}\wedge \dots \wedge \underbrace{I\!I_{f_a}(e_s)}_{j-th} \wedge \dots \wedge e_{i_p})\\
&=&-\sum_{a=1}^m\sum_{s=1}^n\, \sum_{\substack{k,j=1\\ k\neq j}}^p(-1)^{k+1} g(I\!I_{f_a}(e_s),e_{i_j})\delta_{si_k} e_{i_1}\wedge \cdots \wedge{\widehat e_{i_k}} \wedge\cdots \wedge I\!I_{f_a}(X) \wedge \dots \wedge e_{i_p}\\
&&-\sum_{a=1}^m\sum_{s=1}^n\sum_{j=1}^p(-1)^{j+1}  g(I\!I_{f_a}(e_s),e_{i_j})g(I\!I_{f_a}(X),e_s) e_{i_1}\wedge \cdots \wedge{\widehat e_{i_j}} \wedge\cdots \wedge e_{i_p}\\
&&+\sum_{a=1}^m\sum_{s=1}^n \, \sum_{\substack{k,j=1\\ k\neq j}}^p(-1)^{k+1}  g(I\!I_{f_a}(X),e_{i_j})\delta_{si_k} e_{i_1}\wedge \cdots \wedge{\widehat e_{i_k}} \wedge\cdots\wedge I\!I_{f_a}(e_s) \wedge \dots \wedge e_{i_p}\\
&&+\sum_{a=1}^m\sum_{s=1}^n\sum_{j=1}^p (-1)^{j+1}  g(I\!I_{f_a}(X),e_{i_j})g(I\!I_{f_a}(e_s),e_s) e_{i_1}\wedge \cdots \wedge{\widehat e_{i_j}} \wedge\cdots \wedge e_{i_p}.
\end{eqnarray*}
By using the symmetry of the second fundamental form, the above computation reduces to
\begin{eqnarray*}
\sum_{s=1}^n e_s \lrcorner R(e_s,X) \omega&=&-\sum_{a=1}^m \, \sum_{\substack{k,j=1\\ k\neq j}}^p(-1)^{k+1} g(I\!I_{f_a}(e_{i_k}),e_{i_j}) e_{i_1}\wedge \cdots \wedge{\widehat e_{i_k}} \wedge\cdots \wedge I\!I_{f_a}(X) \wedge \dots \wedge e_{i_p}\\
&&-\sum_{a=1}^m\sum_{j=1}^p(-1)^{j+1}  g(I\!I_{f_a}^2(X),e_{i_j}) e_{i_1}\wedge \cdots \wedge{\widehat e_{i_j}} \wedge\cdots \wedge e_{i_p}\\
&&+\sum_{a=1}^m\, \sum_{\substack{k,j=1\\ k\neq j}}^p (-1)^{k+1}  g(I\!I_{f_a}(X),e_{i_j}) e_{i_1}\wedge \cdots \wedge{\widehat e_{i_k}} \wedge\cdots\wedge I\!I_{f_a}(e_{i_k}) \wedge \dots \wedge e_{i_p}\\
&&+\sum_{a=1}^m\sum_{j=1}^p (-1)^{j+1}  g(I\!I_{n\widetilde{H}}(X),e_{i_j}) e_{i_1}\wedge \cdots \wedge{\widehat e_{i_j}} \wedge\cdots \wedge e_{i_p}.
\end{eqnarray*}
Now, we prove that the first sum vanishes. Indeed, by decomposing with respect to $k<j$ and $k>j$, it is equal to
\begin{eqnarray*}
&&\sum_{k < j} (-1)^{k+1} g(I\!I_{f_a}(e_{i_k}),e_{i_j}) e_{i_1}\wedge \cdots \wedge{\widehat e_{i_k}} \wedge\cdots \wedge I\!I_{f_a}(X) \wedge \dots \wedge e_{i_p}\\
&& +\sum_{k > j} (-1)^{k+1} g(I\!I_{f_a}(e_{i_k}),e_{i_j}) e_{i_1}\wedge \cdots \wedge I\!I_{f_a}(X) \wedge\cdots \wedge{\widehat e_{i_k}}  \wedge \dots \wedge e_{i_p}\\
&=&\sum_{k < j} (-1)^{k+j-1} g(I\!I_{f_a}(e_{i_k}),e_{i_j}) I\!I_{f_a}(X) \wedge e_{i_1}\wedge \cdots \wedge{\widehat e_{i_k}} \wedge\cdots \wedge {\widehat e_{i_j}}  \wedge \dots \wedge e_{i_p}\\
&& +\sum_{k >j} (-1)^{k+j} g(I\!I_{f_a}(e_{i_k}),e_{i_j}) I\!I_{f_a}(X) \wedge e_{i_1}\wedge \cdots \wedge {\widehat e_{i_j}} \wedge\cdots \wedge{\widehat e_{i_k}}  \wedge \dots \wedge e_{i_p}\\
&=& 0.
\end{eqnarray*}
In the same way, we prove that the third sum is equal to  $-\sum_{a=1}^mI\!I_{f_a}^{[p-1]}(I\!I_{f_a}(X)\lrcorner\omega)$. Indeed, it is equal to
\begin{eqnarray*}
&&\sum_{k < j} (-1)^{k+1}g(I\!I_{f_a}(X),e_{i_j})(e_{i_1}\wedge \cdots\wedge \widehat{e_{i_k}} \wedge\cdots \wedge I\!I_{f_a}(e_{i_k})\wedge \cdots \wedge e_{i_p})\\
&&+\sum_{k > j} (-1)^{k+1} g(I\!I_{f_a}(X),e_{i_j})(e_{i_1}\wedge \cdots\wedge I\!I_{f_a}(e_{i_k})\wedge \cdots \wedge \widehat{e_{i_k}}\wedge\cdots \wedge e_{i_p})\\
&=&\sum_{k < j} (-1)^{k+1+j-k-1}g(I\!I_{f_a}(X),e_{i_j})(e_{i_1}\wedge \cdots\wedge I\!I_{f_a}(e_{i_k}) \wedge\cdots \wedge {\widehat e_{i_j}} \wedge \cdots \wedge e_{i_p})\\
&&+\sum_{k>j} (-1)^{k+1+k-j-1} g(I\!I_{f_a}(X),e_{i_j})(e_{i_1}\wedge \cdots\wedge \widehat{e_{i_j}} \wedge \cdots \wedge I\!I_{f_a}(e_{i_k})\wedge\cdots \wedge e_{i_p})\\
&=&\sum_{k < j} (-1)^{j}g(I\!I_{f_a}(X),e_{i_j})(e_{i_1}\wedge \cdots\wedge I\!I_{f_a}(e_{i_k}) \wedge\cdots \wedge {\widehat e_{i_j}} \wedge \cdots \wedge e_{i_p})\\
&&+\sum_{k >j} (-1)^{j} g(I\!I_{f_a}(X),e_{i_j})(e_{i_1}\wedge \cdots\wedge \widehat{e_{i_j}} \wedge \cdots \wedge I\!I_{f_a}(e_{i_k})\wedge\cdots \wedge e_{i_p})\\
&=&-I\!I_{f_a}^{[p-1]}(I\!I_{f_a}(X)\lrcorner\omega).
\end{eqnarray*}
In the last equality, we use the formula $A^{[p]}(X_1\wedge \cdots\wedge X_p)=\sum_{i=1}^p X_1\wedge \cdots \wedge A(X_i)\wedge\cdots\wedge X_p$ for any vector fields $X_1,\cdots X_p$ in $TM$.
Therefore, we deduce  that
\begin{eqnarray*}
\sum_{s=1}^n e_s \lrcorner R(e_s,X) \omega&=&-\sum_{a=1}^m I\!I_{f_a}^2(X)\lrcorner \omega-\sum_{a=1}^m I\!I_{f_a}^{[p-1]}(I\!I_{f_a}(X)\lrcorner\omega)+I\!I_{n\widetilde H}(X)\lrcorner\omega\\
&=&- \sum_{a=1}^m I\!I_{f_a}(X)\lrcorner I\!I_{f_a}^{[p]}\omega+I\!I_{n\widetilde H}(X)\lrcorner\omega.
\end{eqnarray*}
In the last equality, we used the identity $A^{[p-1]}(X\lrcorner\alpha)=X\lrcorner A^{[p]}\alpha-A(X)\lrcorner \alpha$. 
\hfill$\square$

\end{document}